\documentclass[11pt,reqno,a4paper,english]{amsart}

\usepackage[utf8]{inputenc}  
\usepackage{dsfont}
\usepackage{xcolor}

\usepackage{enumerate}          
\usepackage{amssymb,amsmath,amsfonts,amsthm}            
\usepackage{url}
\usepackage{hyperref}
\usepackage{soul}
\usepackage{anysize}
\usepackage{mathrsfs}
\usepackage{graphicx,float,color}
\usepackage{epsfig}
\usepackage{tikz}
\usepackage{caption}
\usetikzlibrary{patterns, patterns.meta}
\usetikzlibrary{calc,arrows.meta,positioning}
\usetikzlibrary{decorations.pathreplacing}

\usepackage{hyperref}

\hypersetup{
   colorlinks,
   linkcolor={red},
   citecolor={blue},
   urlcolor={red}
}

\usepackage{macros2}

\usepackage{cite}

\usepackage{marginnote}

\AtBeginDocument{%
  \setlength{\abovedisplayskip}{9pt plus 2pt minus 2pt}%
  \setlength{\belowdisplayskip}{9pt plus 2pt minus 2pt}%
  \setlength{\abovedisplayshortskip}{5pt plus 2pt}%
  \setlength{\belowdisplayshortskip}{6pt plus 2pt minus 2pt}%
}

\begin{document}

\theoremstyle{plain}
\newtheorem{theorem}{Theorem}[section]
\newtheorem{definition}[theorem]{Definition}
\newtheorem{proposition}[theorem]{Proposition}
\newtheorem{lemma}[theorem]{Lemma}

\newtheorem{corollary}[theorem]{Corollary}
\newtheorem{conjecture}[theorem]{Conjecture}
\newtheorem{remark}[theorem]{Remark}
\newtheorem{open}[theorem]{Open question}
\numberwithin{equation}{section}
\newtheorem*{claim}{Claim}
\newtheorem{assumption}[theorem]{Assumption}
\errorcontextlines=0

\newcommand{\Id}{\text{Id}}
\newcommand{\mnum}{A}
\newcommand{\mrat}{B}
\newcommand{\R}{\mathbb{R}}
\newcommand{\Z}{\mathbb{Z}}
\newcommand{\N}{\mathbb{N}}
\newcommand{\X}{\mathcal{X}}
\newcommand{\Y}{\mathcal{Y}}
\newcommand{\BW}{{\rm BW}}
\newcommand{\red}{\color{red}}
\newcommand{\black}{\color{black}}
\newcommand{\A}{A}
\newcommand{\Tr}{{\rm Tr}}
\newcommand{\ppow}{\beta}
\newcommand{\gau}{\ppow}
\newcommand{\Cor}{\mathcal{C}}
\renewcommand{\H}{\mathcal{H}}
\newcommand{\B}{B}
\newcommand{\Ball}{\mathcal{B}}
\newcommand{\C}{C}
\newcommand{\D}{D}
\newcommand{\E}{E}
\newcommand{\ct}{5\B\sqrt{d}}
\newcommand{\pow}{p}
\newcommand{\powun}{p_1}
\newcommand{\powde}{p_2}
\newcommand{\pds}{\delta}
\renewcommand{\P}{\mathbb{P}}
\renewcommand{\S}{\mathcal{S}}
\newcommand{\dist}{{\rm dist}}
\newcommand{\spt}{{\rm spt}}
\renewcommand{\diam}{\operatorname{diam}}
\newcommand{\inte}{{\rm int}}
\newcommand{\Per}{{\rm Per}}
\newcommand{\vol}{{\rm vol}}
\newcommand{\CDM}{C_{\rm DM}}
\newcommand{\Res}{\mathcal{B}}
\newcommand{\Kant}{\mathcal K}
\newcommand{\YSp}{\mathcal Y}
\newcommand{\Rsp}{\mathbb R}
\newcommand{\eps}{\varepsilon}
\newcommand{\Class}{\mathcal C}
\newcommand{\Conv}{\mathrm{Conv}}
\newcommand{\sca}[2]{\langle #1\vert #2\rangle}
\newcommand{\nr}[1]{\Vert #1\Vert}
\newcommand{\restr}[1]{\left. #1\right\vert}
\newcommand{\dd}{{\rm d}}
\renewcommand{\Box}{{\rm Box}}

\hyphenation{mar-gi-nals}

\setcounter{tocdepth}{1}

\title[Sharp Stability of Brenier Potentials and its Link to Poincaré Inequalities]{Sharp Stability of Brenier potentials and \\its link to Poincaré inequalities}

\author{Cyril Letrouit}\address{Cyril Letrouit. Université Paris-Saclay, CNRS, Inria, Laboratoire de mathématiques d’Orsay, 91405, Orsay, France} \email{cyril.letrouit@universite-paris-saclay.fr}

\author{Tudor Manole}\address{Tudor Manole. Stanford University, Department of Statistics, California, United States} \email{tmanole@stanford.edu}

\date{\today}

\maketitle

\begin{abstract}
We prove a new functional inequality for uniformly log-concave measures, stronger than the Poincaré inequality. This new inequality bounds the variance of differences of convex functions and is motivated by the stability theory of optimal transport: for a fixed uniformly log-concave source measure,
it provides a sharp stability bound for Brenier potentials with respect to variations of the target measure. 
Applying this inequality to Gaussian measures, we recover the fact that the Bures-Wasserstein distance between two positive semidefinite matrices is equivalent to the Frobenius distance between their square roots. 
Our proofs hinge upon linearization arguments
for the Monge-Amp\`ere equation governing
the quadratic optimal transport problem.
\end{abstract}

{
    \hypersetup{hidelinks}
    \tableofcontents
}

\section{Introduction}
Functional inequalities, 
like the Poincar\'e, logarithmic
Sobolev, and transportation
inequalities, play
a fundamental role in high-dimensional
probability theory and convex geometry.
Among others, they provide important  tools
for the study of concentration of measure
and mixing times of Markov processes~\cite{ledoux2001concentration,bakry2014analysis},
with applications 
to high-dimensional sampling, Bayesian statistics,
statistical physics, and machine learning \cite{villani2003,Chewi24}. 
The goal of this manuscript
is to introduce a new functional inequality 
motivated
by the recent
literature on quantitative
stability of optimal transport
maps, and to discuss
its connections to Poincar\'e
inequalities.

Recall that a probability
measure $\rho$ on $\bbR^d$, absolutely continuous with respect to the Lebesgue measure,
is said to satisfy 
a Poincar\'e inequality
if there exists a finite 
constant $C > 0$ such that
\begin{equation}\label{e:poincarerho}
\forall f \in C_c^\infty(\R^d), \qquad \Var_\rho(f)\leq C \|\nabla f\|^2_{L^2(\rho)}
\end{equation}
where $\Var_\rho(f)=\int f^2\,\dd\rho-(\int f\,\dd\rho)^2$ is the variance with respect to $\rho$, and $C_c^\infty(\R^d)$ denotes the set of smooth functions with compact support. We denote by  $C_{P}(\rho)$ the
Poincar\'e constant of $\rho$, 
namely the 
smallest constant $C$
for which inequality~\eqref{e:poincarerho}
holds. 

Log-concave probability
distributions form
an important class of measures
that satisfy Poincar\'e inequalities. 
Recall that $\rho $
is said to be log-concave 
if it admits  a density 
with respect to the Lebesgue
measure on $\bbR^d$
which is 
proportional to $e^{-V}$, where $V$ is a convex function on $\bbR^d$,
and is  further 
said to be $\kappa$-uniformly log-concave  if $V-\frac{\kappa}{2}|\cdot|^2$ is convex for some $\kappa > 0$.
A classical result
asserts that 
any $\kappa$-uniformly
log-concave measure satisfies
the Poincar\'e inequality
\begin{equation}\label{eq:unif_Cp}
C_P(\rho) \leq \frac 1 \kappa.
\end{equation}
More generally,
recent breakthroughs
on the Kannan-Lov\'asz-Simonovits conjecture (cf.~\cite{bizeul2026presenting}
and references therein)
imply that the Poincaré constant $C_P(\rho)$ of any log-concave 
measure $\rho$ is upper bounded by the operator norm of the covariance of $\rho$, multiplied by a universal constant.
 
Let $\mathcal{P}_2(\R^d)$ denote the set of probability measures with finite second moment. The main object of interest in this
paper is the following 
functional inequality: We shall
say that an absolutely continuous probability measure
$\rho\in\mathcal{P}_2(\R^d)$ satisfies a
{\it Brenier stability inequality} if there
exists a finite constant $C > 0$ such that 
\begin{equation}
\label{eq:stability}
\begin{aligned}
&\Var_\rho(\phi_1-\phi_0)\leq C W_2^2((\nabla\phi_0)_\#\rho, (\nabla\phi_1)_\#\rho),
\end{aligned}
\end{equation}
for all 
convex functions 
$\phi_0,\phi_1$ on $\bbR^d$
satisfying 
$\nabla\phi_0,\nabla\phi_1 \in L^2(\rho)$. Here, $W_p$ denotes the 
$p$-Wasserstein distance, 
and $F_\# \rho$ denotes the pushforward
of $\rho$ by a measurable
vector field $F$. 
Further notational conventions
are recalled in~Section~\ref{sec:notation}.
We denote by $C_W(\rho)$
the smallest
constant $C > 0$ 
for which equation~\eqref{eq:stability}
holds. If no such $C<\infty$ exists, we write $C_W(\rho)=\infty$.

The Brenier stability inequality is stronger than the Poincaré inequality: 
for any $\rho \in \calP_2(\bbR^d)$,
\begin{equation}\label{eq:brenier_implies_poincare}
C_P(\rho) \leq C_W(\rho).
\end{equation}
Indeed, since $(\nabla\phi_0,\nabla\phi_1)_\#\rho$ is a coupling between $(\nabla\phi_0)_\#\rho$ and $(\nabla\phi_1)_\#\rho$, the following well-know inequality holds: \begin{equation}\label{e:coupling}
W_2((\nabla\phi_0)_\#\rho,(\nabla\phi_1)_\#\rho)\leq \|\nabla\phi_0-\nabla\phi_1\|_{L^2(\rho)}.
\end{equation}
Hence 
equation~\eqref{eq:stability} implies equation~\eqref{e:poincarerho} with same constant $C$, for any $f$ of the form $f=\phi_1-\phi_0$, where $\phi_0,\phi_1$ are convex and satisfy 
$\nabla\phi_0,\nabla\phi_1\in L^2(\rho)$. To deduce
the Poincar\'e inequality~\eqref{e:poincarerho}  for arbitrary
$f \in \calC_c^\infty(\bbR^d)$, we can apply the same
argument to $\phi_0=M|\cdot|^2$ and $\phi_1=f+M|\cdot|^2$,  for a constant $M\geq 0$ 
chosen large enough so that $\phi_1$ is convex. Then, both $\phi_0$ and $\phi_1$
are convex 
and satisfy $\nabla\phi_0,\nabla\phi_1\in L^2(\rho)$, hence \eqref{e:poincarerho} holds with $C=C_W(\rho)$, which proves \eqref{eq:brenier_implies_poincare}. 

Another way to see the relationship between the Brenier stability inequality and Poincaré-type inequalities is to use the following known characterization of the 2-Wasserstein distance, whose proof is recalled in Appendix \ref{a:caracW2} (see also~\cite[Lemma 2.7]{pinzi2025totally}): For all convex maps $\phi_1,\phi_0$
with gradients in $L^2(\rho)$, one has  
\begin{equation}\label{e:W2asinf}
W_2((\nabla\phi_0)_\#\rho , (\nabla\phi_1)_\#\rho) = \inf_{S \in \calS(\rho)} \|\nabla \phi_1 -   \nabla \phi_0\circ S \|_{L^2(\rho)}
\end{equation}
where $\calS(\rho)$ denotes the group of measure-preserving maps
of $\rho$, i.e., the Borel maps $S:\R^d\rightarrow \R^d$ such that $S_\#\rho=\rho$. As a consequence of this identity, a measure $\rho \in \calP_2(\bbR^d)$
satisfies a Brenier stability inequality
if and only if it satisfies the bound
\begin{align}
\label{eq:stability_polar}
\Var_\rho(\phi_1-\phi_0) \leq C_W(\rho) \cdot \|\nabla \phi_1-  \nabla \phi_0\circ S\|_{L^2(\rho)}^2,
\end{align} 
for all $S \in \calS(\rho)$ 
and all convex functions $\phi_1,\phi_0$
whose gradients lie in $L^2(\rho)$. Again, for $S={\rm id}$, we recover the Poincaré inequality and \eqref{eq:brenier_implies_poincare}, but for  general $S\in \calS(\rho)$, equation~\eqref{eq:stability_polar} provides a more general Poincaré-type inequality.

\subsection{Main Results}
The purpose of this manuscript is to 
prove that the Brenier stability inequality 
holds under natural assumptions on $\rho$. We will additionally show that for some (non-uniformly) log-concave $\rho$, it can happen that $C_P(\rho)<C_W(\rho)=\infty$, thus showing that the Brenier stability inequality is \emph{strictly} stronger than the Poincaré inequality. Our first main result is the following.
\begin{theorem}\label{t:main}
Let $\rho$ be a $\kappa$-uniformly log-concave probability measure on $\R^d$. Then, for any convex
functions $\phi_0,\phi_1$
on $\bbR^d$ such that $\nabla \phi_0,\nabla\phi_1\in L^2(\rho)$, there holds 
\begin{equation}\label{e:main}
\Var_\rho(\phi_1-\phi_0)\leq \frac 1 \kappa W_2^2((\nabla\phi_0)_\#\rho, (\nabla\phi_1)_\#\rho).
\end{equation}
That is, $C_W(\rho) \leq 1/\kappa$. 
\end{theorem}
In other words, any $\kappa$-uniformly log-concave probability measure $\rho$ satisfies not only $C_P(\rho)\leq1/\kappa$, but even $C_W(\rho)\leq 1/\kappa$.
It is a classical fact that $C_P(\rho)=1$ for standard Gaussians
\cite{chernoff1981note,bakry2014analysis}, 
and in this case, Theorem~\ref{t:main}
together with the inequality~\eqref{eq:brenier_implies_poincare} imply  the remarkable identity 
\begin{equation}
\label{eq:gaussian_constant}
C_P(\rho)=C_W(\rho)=1.
\end{equation}
This equality
has notable
consequences for the Bures-Wasserstein
geometry, which we
describe in Section~\ref{sec:intro_bures} below.

Theorem~\ref{t:main} raises
the question
of whether all log-concave
measures 
satisfy a Brenier stability inequality. 
This turns out not be the case, as we
show by using the following 
characterization of Brenier stability inequalities
for unconditional measures. Recall that a measure $\rho$ in $\R^d$ is said to be unconditional if it is invariant under coordinate reflections, i.e., for
any $g\in C_c^\infty(\R^d)$ and any choice of signs,
$$ \int g(\pm x_1,\pm x_2, \ldots, \pm x_d) \dd\rho(x) = \int g(x_1,\ldots,x_d)\dd\rho(x).$$ 
\begin{theorem}\label{t:directionalpoincare}
Suppose that $\rho\in \calP_2(\bbR^d)$
is an unconditional and absolutely continuous
probability measure satisfying
$C_W(\rho)< \infty$. 
Suppose further that the density
of $\rho$ is twice continuously
differentiable on $\bbR^d$, and strictly
positive.
Then, for any $j=1,\dots,d$,  
$\rho$ satisfies the directional Poincar\'e
inequality
\begin{align} \label{eq:directional_poincare}
\Var_\rho(f)\leq C_W(\rho) \int \left|\frac{\partial f}{\partial x_j} (x)\right|^2d\rho(x),
\end{align}
for any function 
$f\in \calC_c^\infty(\bbR^d)$ that is odd in its $j$-th variable.
\end{theorem}
Theorem~\ref{t:directionalpoincare}
is proven using a coupling
argument in Section~\ref{sec:negative}.
It turns out that there exist unconditional
log-concave measures which fail to satisfy
the directional Poincar\'e inequality~\eqref{eq:directional_poincare}, 
thus such measures  cannot satisfy
a Brenier stability inequality. We summarize
this conclusion as follows.
\begin{corollary}\label{c:failkls}
For any $d \geq 2$, 
there exists  a 
twice differentiable convex function 
$V:\bbR^d \to \bbR$
satisfying for some $\lambda > 0$  
$$\forall x\in\R^d, \qquad \nabla^2 V(x) \preceq \lambda I_d$$
such that the log-concave measure $\rho$ with density $ e^{-V}$ with respect to the Lebesgue measure
satisfies $C_W(\rho)=\infty.$
\end{corollary}
Since all log-concave
probability measures have a finite Poincar\'e
constant, Corollary~\ref{c:failkls}
shows that there are examples of measures
which do not satisfy a Brenier stability inequality 
but do satisfy an $L^2$ Poincar\'e inequality. 

Our original motivation for proving Theorem \ref{t:main} comes from the quantitative stability theory of optimal transport. We briefly recall the aims of this area here and defer a more detailed discussion to Section \ref{s:stabbrenier}.  
Given $\rho,\mu\in\mathcal{P}_2(\R^d)$ such that $\rho$ is absolutely continuous with respect to the Lebesgue measure, Brenier's theorem \cite{brenier1991} asserts that there exists a $\rho$-almost everywhere unique solution $T:\bbR^d\to\bbR^d$ to the optimization problem
\begin{equation}\label{e:ot}
\inf_{T \, : \, T_\#\rho = \mu} \int |T(x)-x|^2 \dd\rho,
\end{equation}
which is  nowadays called the Brenier map.
Moreover, the map $T$ takes the form $T=\nabla\phi$ for some convex function $\phi:\R^d\rightarrow \R$ called a Brenier potential. In the rest of this discussion, 
we treat the source measure 
$\rho$ as fixed, and we assume that for any $\mu$, 
the Brenier potential is unique up to addition of a constant, which is the case for instance when $\rho$ satisfies the Poincar\'e inequality \eqref{e:poincarerho} (as a consequence of the fact that $T=\nabla\phi$ itself is $\rho$-a.e. unique). 
To emphasize the dependence on the target measure $\mu$, we denote by $T_{\mu}$ the  solution to 
the optimization problem~\eqref{e:ot}, and by $\phi_{\mu}$ the unique Brenier potential satisfying $\int \phi\, \dd\rho =0$. 

It is natural to ask whether the map $T_\mu$ or the potential
$\phi_\mu$ are {\it stable} with respect to perturbations of
the target measure $\mu$. 
Qualitative answers to this question have
been well-understood
since Brenier's original paper \cite{brenier1991}, where it
was shown that the map $\mu\mapsto T_\mu$ is continuous from $\mathcal{P}_2(\R^d)$ endowed with the $2$-Wasserstein distance $W_2$ to $L^2(\rho)$. When $\rho$ satisfies the Poincaré inequality, this implies continuity of $\mu\mapsto \phi_\mu$ as well.
Quantitative versions of these stability results, by contrast, have attracted attention only more recently, as reviewed in Section~\ref{s:stabbrenier}.
For Brenier potentials, existing quantitative results typically take the form
\begin{equation}\label{e:linear}
\forall \mu,\nu\in\mathcal{P}(Y), \qquad \|\phi_\mu-\phi_\nu\|_{L^2(\rho)}\leq CW_p(\mu,\nu)^\alpha
\end{equation}
where $p\geq 1$, $\mathcal{P}(Y)$ denotes the set of probability measures supported in a given compact set $Y$, and  $C,\alpha>0$ are constants which may depend on $Y$, but not on $\mu,\nu\in \mathcal{P}(Y)$. 
In this context, Theorem \ref{t:main} shows that if $\rho$ is uniformly log-concave, then 
equation~\eqref{e:linear} holds for $\alpha=1$ and $p=2$, and $\mathcal{P}(Y)$ can even be replaced by the larger set $\mathcal{P}_2(\R^d)$.
On the other hand, Corollary~\ref{c:failkls} shows that the Lipschitz exponent
$\alpha=1$ is generally unachievable in equation~\eqref{e:linear} for log-concave source measures.

It was conjectured in \cite[Conjecture 3.10]{letrouit2025quantitative} that if $\rho$ is uniform on a convex body (i.e., a compact, convex set, with non-empty interior), then for any bounded set $Y$, equation~\eqref{e:linear} holds for $\alpha=1$, $p=2$, and some $C>0$. A proof of this conjecture
was recently claimed in OpenAI's document \cite[Proposition 5.1]{openai}, which does not mention \cite[Conjecture 3.10]{letrouit2025quantitative}. We provide in this paper a different proof, close to the one of Theorem \ref{t:main}.
\begin{theorem}\label{t:convexset}
Let $K\subset B(x_0,R)$ be a convex body, and let
$$
\dd\rho=e^{-V(x)}\mathbf{1}_K(x)\,\dd x,
$$
where $\mathbf{1}_K$ is the characteristic function of $K$, and $V$ is convex. Then for any  convex functions $\phi_0,\phi_1$
on $\bbR^d$ such that $\nabla\phi_0,\nabla\phi_1\in L^2(\rho)$, there holds
\begin{equation}\label{e:mainconvex}
\Var_\rho(\phi_0-\phi_1)
\leq
4R^2\,W_2^2((\nabla\phi_0)_\#\rho,(\nabla\phi_1)_\#\rho).
\end{equation}
That is, $C_W(\rho)\leq 4R^2$.
\end{theorem}
As a corollary of Theorem \ref{t:convexset}, we obtain the following generalization of the sharp $1/3$-stability result for optimal transport maps claimed in OpenAI's document \cite[Theorem 1.1]{openai}.
\begin{corollary}\label{c:1/3}
Let $K\subset B(x_0,R)$ be a convex body, and let
$$
\dd\rho=e^{-V(x)}\mathbf{1}_K(x)\,\dd x,
$$
where $V$ is convex and takes finite values on $K$. Then for any compact set $Y\subset\R^d$, there exists $C>0$ such that
\begin{equation}\label{e:mainconvexmap}
\forall \mu_0,\mu_1\in\mathcal{P}(Y), \qquad \|T_{\mu_0}-T_{\mu_1}\|_{L^2(\rho)}
\leq
C\,W_2(\mu_0,\mu_1)^{1/3}.
\end{equation}
Here $\mathcal{P}(Y)$ denotes the set of probability measures on $Y$, and, for $i\in\{0,1\}$, $T_{\mu_i}$  denotes the quadratic optimal transport map from $\rho$ to $\mu_i$.
\end{corollary}
Corollary~\ref{c:1/3} is a consequence of
Theorem~\ref{t:convexset}, 
together with the reverse Poincar\'e inequality
established in~\cite[Section 4]{delalande2023quantitative}.
We present a short proof in Section~\ref{s:coercive_cvx}.

\subsection{Consequences for the Bures-Wasserstein Geometry}
\label{sec:intro_bures}
In the particular case where~$\rho$ is Gaussian and $\phi_0,\phi_1$ are quadratic, Theorem \ref{t:main} has an interesting corollary.
Recall that the Bures–Wasserstein space $\BW(\R^d)$ is the space of non-degenerate Gaussians on $\R^d$, equipped with the Wasserstein metric. Over this space, 
the squared 
Wasserstein distance takes on a closed form, given by
\begin{equation}\label{e:BWformula}
W_2^2(\mathcal{N}(m_0,\Sigma_0),\mathcal{N}(m_1,\Sigma_1))= |m_0-m_1|^2 + B^2(\Sigma_0,\Sigma_1),
\end{equation}
where $\mathcal{N}(m,\Sigma)$ denotes the Gaussian with mean $m$ and covariance $\Sigma$, and
 $B(\Sigma_0,\Sigma_1)$
denotes the Bures metric over the space
of positive semidefinite matrices, defined by
\begin{equation}
B(\Sigma_0,\Sigma_1) = \sqrt{\Tr(\Sigma_0)+\Tr(\Sigma_1)-2\Tr[(\Sigma_0^{1/2}\Sigma_1\Sigma_0^{1/2})^{1/2}]}.
\end{equation}
The Bures-Wasserstein geometry has been the subject of many works; see for instance \cite{bhatia2019bures, bhatia2019matrix}. The following corollary shows that the  
Bures distance
is comparable up to a factor $2$ to the simpler
metric
$$H(\Sigma_0,\Sigma_1)
:=\|\Sigma_0^{1/2}-\Sigma_1^{1/2} \|_{F},$$
where $\|\cdot\|_F$ denotes the Frobenius norm. It was already proved by other means in \cite[Proposition 2]{dinh2020some} and \cite[Proposition 2]{dinh2021weighted}, 
and also has a known
counterpart in the quantum information theory literature~\cite[Proposition 9]{spehner2017geometric}, 
but does not appear to be 
widely-known 
in the optimal transport
literature.
\begin{corollary} \label{c:bures} For any $m_0,m_1\in\R^d$ and $\Sigma_0,\Sigma_1$ positive definite matrices, there holds
$$
\frac12 H^2(\Sigma_0,\Sigma_1) \leq W_2^2(\mathcal{N}(m_0,\Sigma_0),\mathcal{N}(m_1,\Sigma_1))
-|m_0-m_1|^2\leq H^2(\Sigma_0,\Sigma_1).
$$
\end{corollary}
This result is proved in   Section \ref{s:bures} as a direct corollary of Theorem \ref{t:main}, 
using the observation that if $\phi_0,\phi_1$
are Brenier potentials
whose gradients respectively push forward the standard
Gaussian measure $\rho=\calN(0,I_d)$ onto
$\calN(m_0,\Sigma_0)$ and $\calN(m_1,\Sigma_1)$,
then 
$$\Var_\rho(\phi_1-\phi_0) = \frac12 H^2(\Sigma_0,\Sigma_1) + |m_0-m_1|^2.$$
Corollary \ref{c:bures} then follows
using the fact that standard Gaussians satisfy
a Brenier stability inequality
with constant $C_W(\rho)=1$, as stated in equation~\eqref{eq:gaussian_constant}.

\subsection{Connections to Stability
of Optimal Transport} \label{s:stabbrenier} 
As already mentioned, one of the main applications of our results is to the quantitative stability theory of optimal transport. More precisely, Theorem \ref{t:main} shows that if $\rho$ is $\kappa$-uniformly log-concave, then
Brenier potentials are linearly stable with respect to $W_2$.
Said differently, the map taking $\mu$ to $\phi_\mu$ defines
a Lipschitz embedding of the Wasserstein space $(\calP_2(\bbR^d),W_2)$
into $L^2(\rho)$.
In contrast, 
Corollary~\ref{c:failkls} shows that this is no longer the case if $\rho$ is merely assumed to be log-concave.

As Corollary \ref{c:1/3} already suggests, the quantitative stability of Brenier potentials is closely related to the quantitative stability of Brenier maps, i.e., to bounds of the form 
\begin{equation}\label{e:stabmap}
\forall \mu,\nu\in\mathcal{P}(Y), \qquad \|\nabla \phi_\mu-\nabla\phi_\nu\|_{L^2(\rho)}\leq CW_p(\mu,\nu)^\beta
\end{equation}
for $C,\beta>0$, $p\geq 1$ and $Y\subset\R^d$.
Both problems have been the subject of intensive research in recent years. 
In what follows, we briefly review this literature through the lens of Brenier potentials, since these are the objects investigated in the present paper. We refer to~\cite{letrouit2025quantitative}
for a broader survey.

There are two predominant sets of assumptions under
which such stability bounds have been analyzed. 
First, strong results have been achieved under regularity assumptions on one of the two Brenier maps. In this line of work, one establishes \eqref{e:stabmap} under the assumption that $\nabla\phi_\mu$ is Lipschitz (and \eqref{e:linear} then follows if $\rho$ satisfies the Poincaré inequality). The constant $C$ in \eqref{e:stabmap} then depends on the Lipschitz constant of $\nabla\phi_\mu$.
The first result of this type was proved in Gigli's work \cite{gigli2011holder},
who obtained an exponent of $\beta=1/2$, and several variants of this result
were later derived in the literature on statistical optimal transport~\cite{hutter2021,ghosal2022multivariate,deb2021rates,vacher2024optimal}. 
When $\phi_\mu$ is additionally assumed to be strongly convex, it is 
known that the exponent $\beta=1$ can be achieved in equation~\eqref{e:stabmap},
cf.~\cite{manole2024plugin,balakrishnan2026stability,gonzalez2026empirical}.

A second
line of work  establishes bounds of the form \eqref{e:linear} or \eqref{e:stabmap} under the only assumption that $Y$ is compact. No other hypothesis is made on $\mu$ and $\nu$, or their induced
potentials $\phi_\mu$ and $\phi_\nu$. We can mention the works \cite{berman2021convergence,delalande2023quantitative,letrouitmerigot} which show that 
if~$Y$ is compact and $\rho$ is either uniformly log-concave or has a density bounded above and below on a nice domain (typically, a John domain), then \eqref{e:linear} holds 
with $p=1$ and $\alpha=1/2$ (up to a multiplicative logarithmic term if $\rho$ is uniformly log-concave). Notably, the article \cite{letrouitmerigot} was the first to point out the relation between stability of Brenier potentials and Poincaré-type inequalities. For $p=1$, the exponent $\alpha=1/2$ is optimal as shown in \cite[Proposition 6.10]{letrouitmerigot} (see also \cite[Section 3.7]{letrouit2025quantitative}). Of course, since $W_1\leq W_2$, these works imply the same bounds for $p=2$ instead of $p=1$ in \eqref{e:linear}. However, for $p=2$, the exponent $\alpha=1/2$ is no longer optimal. The present work and OpenAI's document \cite{openai} are the first   to establish the optimal exponent when $p=2$, namely $\alpha=1$. Notice however that \cite{openai} does not handle uniformly log-concave measures and provides a worse constant than $4R^2$ in Theorem \ref{t:convexset}.

The present paper falls under the second line of work, since  we make no assumptions on~$\mu,\nu$. So far, the predominant proof strategy of the second line of works was to establish some form of strong convexity for the  Kantorovich functional defined by $$ \mathcal{K}_\rho:\psi\mapsto \int \psi^* \, \dd\rho, $$ with the notable exception of \cite{divol2025tight} which built on stability of entropic optimal transport maps. The recent document of OpenAI~\cite{openai} exploits another approach,  based on semi-discrete computations and displacement concavity. The main results of the present paper are proved using fundamentally different arguments, inspired from the first line of work, notably from \cite{manole2023,gonzalez2026empirical}. Indeed, our proofs hinge upon the perturbation theory for the Monge-Amp\`ere equation, which until now has typically relied on some form of regularity on $\mu$ and $\nu$.  Our proof strategy is summarized in Section \ref{s:proofstrategy}. 

\subsection{Outline}
The remainder of this manuscript
is organized as follows. In Section \ref{s:proofstrategy}, we present an overview of the proofs. In Section \ref{s:reductocoercest}, we introduce the linearized Monge-Ampère operator $L_\mu$ and show how a coercivity estimate on $L_\mu$, denoted by $A(\kappa)$, leads to our main results. In Section \ref{s:reduccoerciveEphi}, we show that $A(\kappa)$ is implied by a coercivity estimate $B(\kappa)$ on an operator $E_\mu$ obtained by pulling back $L_\mu$ to the source domain, and which has the advantage of being symmetric. We thus reduce to proving $B(\kappa)$.
In Section~\ref{s:coercive_ulc},
we establish Theorem~\ref{t:main} (and deduce Corollary \ref{c:bures}) by showing that $B(\kappa)$ holds for any uniformly log-concave measure $\rho$. This relies on a generalization of Bochner's (and Reilly's) formula. In Section \ref{s:coercive_cvx}, we establish 
Theorem~\ref{t:convexset} and Corollary \ref{c:1/3} by showing that $B(\kappa)$ holds for any compactly supported log-concave $\rho$. This coercivity estimate is actually slightly weaker than for \emph{uniformly} log-concave measures, and depends on the size of the support of $\rho$.
In Section~\ref{sec:negative}, we then prove
Theorem~\ref{t:directionalpoincare} and Corollary~\ref{c:failkls},
establishing a counter-example to
the Brenier stability inequality for arbitrary log-concave probability
measures. In Appendix \ref{a:caracW2}, we prove the formula \eqref{e:W2asinf} for the Wasserstein distance. In Appendix \ref{app:regrhoandK}, we provide a regularization argument allowing to reduce the proof of our main results to the case where $V$ is smooth and $K$ has smooth, uniformly convex boundary.

\subsection{Notation}
\label{sec:notation}
We denote by $|\cdot|$ the Euclidean norm, and by $\langle  \cdot,\cdot\rangle$ the Euclidean scalar product. We let $I_d$ be the identity matrix, and write $A\succeq B$ for two real
symmetric matrices when $A-B$ has
nonnegative eigenvalues. Given a convex function $\varphi:\bbR^d\to\bbR$, 
we let 
\begin{equation}\label{e:convexconj}
\varphi^*(x) = \sup_{y\in \bbR^d} \{ \langle x,y\rangle-\varphi(y)\}
\end{equation}
denote its convex conjugate.

Let $\calP_p(\bbR^d)$ denote
the set of Borel probability
measures on $\bbR^d$ with finite
$p$-th moment for any given
$p\geq1$. For any measures $\mu,\nu\in\mathcal{P}_p(\R^d)$, 
 define
the $p$-Wasserstein distance
$$W_p^p(\mu,\nu) = \inf_{\pi\in\Pi(\mu,\nu)}
\int |x-y|^p \dd\pi(x,y),$$
where $\Pi(\mu,\nu)$ denotes
the set of all couplings of $\mu$ and $\nu$.

Let $\calC^s(A)$ be the standard
H\"older spaces on a subset $A\subseteq \bbR^d$ for integers
$s \geq 1$.
Furthermore, let $\calC_c^\infty(\bbR^d)$ denote
the set of   compactly-supported
infinitely-differentiable
functions on $\bbR^d$.

We recall that under mild assumptions which are satisfied by all probability measures $\rho$ considered in the present paper, the weighted Sobolev space $W^{1,2}(\rho)$ is defined as follows. Assume that there exists an open set $\Omega$ with $\rho > 0$ a.e. on $\Omega$ and $\rho = 0$ a.e. outside $\Omega$. Assume also that $\rho^{-1}$ is locally integrable in $\Omega$. Under these conditions, any $f\in L^2(\rho)$ is locally integrable in $\Omega$ by Cauchy-Schwarz, thus their distributional derivative makes sense. Then, $W^{1,2}(\rho)$ is defined as the set of $f \in L^2(\rho)$ such that their distributional derivative $\nabla f$ is in $L^2(\rho)$, and we let
$$
\|f\|_{W^{1,2}(\rho)}=\left(\|f\|_{L^2(\rho)}^2+\|\nabla f\|_{L^2(\rho)}^2\right)^{1/2}.
$$

Given a measure $\rho \in \calP_2(\bbR^d)$ admitting a 
continuously differentiable density, we define
the weighted divergence and Laplacian operators 
$$\div_\rho F := \frac 1 \rho\div(\rho F),
\quad \text{and}\quad 
\Delta_\rho u = \div_\rho(\nabla u)
= \Delta u - \langle \nabla V, \nabla u\rangle,  $$
for all  $u \in \calC^2(\bbR^d)$ 
and $F \in \calC^1(\bbR^d;\bbR^d)$.
Finally,  let
$
L_0^2(\rho) = \{u
\in L^2(\rho) \mid  \textstyle\int u \, \dd\rho =0 \}.
$

\subsection*{Statement on AI use.}

Our collaboration started in March 2026 while both authors were visiting the IAS in Princeton. From March to early October 2026, we worked toward
proving C.L.'s conjecture \cite[Conjecture 3.10]{letrouit2025quantitative} and Theorem \ref{t:main} without using   AI tools,  apart from searching for a counterexample to the conjecture.
During that period, we 
established Corollary~\ref{c:bures} via independent means, 
developed our current interpretation  
which links stability of optimal transport
to  Poincaré inequalities, 
and developed two proof strategies.
One of these strategies
was the starting point
for the   proof presented here, 
building upon T.M.'s previous work on the topic~\cite{manole2023}.  
In particular, we 
formulated condition $A(\kappa)$ and understood that it would be
 a key step to proving our main results.
On October 6, 2026, OpenAI released the document \cite{openai}, which claims in its Proposition 5.1 a proof of C.L.'s  \cite[Conjecture 3.10]{letrouit2025quantitative}, based on a different approach via semi-discrete optimal transport. 
After this point, we decided to use AI to 
find the arguments that we were missing to 
complete our proof strategy. Our most significant
use of AI was for establishing condition $A(\kappa)$ under a (uniform) log-concavity assumption on $\rho$. We then simplified, reformulated, and extended the proof.  
The writing of the paper is entirely due to the authors. We take full responsibility for the content of the article.

\subsection*{Acknowledgments.} 
This collaboration began while the authors
were visiting the Institute for
Advanced Study, and they would like
to thank the Institute for its hospitality. 
C.L. acknowledges funding from PEPS JCJC and the support of the Agence nationale de la recherche, through the PEPR PDE-AI project (ANR-23-PEIA-0004).

\section{Proof strategy}\label{s:proofstrategy}

 The goal of this section is to provide a heuristic
 overview of the proof
 of our main result, Theorem~\ref{t:main}. 
We assume for the remainder
of this section  that the source
measure $\rho$ is uniformly log-concave, with smooth density over $\R^d$, namely for some $\kappa>0$,
$$\dd\rho(x) = e^{-V(x)}\dd x, \quad \text{with }
\nabla^2 V(x)\succeq \kappa I_d
~\text{ for any } x\in\R^d.$$ 
The weighted negative Laplacian $-\Delta_\rho$ is 
essentially self-adjoint on $C_c^\infty(\mathbb{R}^d)$ in $L^2(\rho)$, and has  compact resolvent~\cite{bakry2014analysis}. It
admits a real and discrete spectrum, which can be used
to define the fractional Laplacian $(-\Delta_\rho)^{s/2}$
for any $s \in \R$.
We can then 
define the 
spectral homogeneous
Sobolev norms 
$$\forall u\in \calC_c^\infty(\bbR^d)\cap L_0^2(\rho), 
\quad \|u\|_{\dot  H^s(\rho)} 
:= \big\|(-\Delta_\rho )^{s/2} u\big\|_{L^2(\rho)},$$
where 
 we recall that $L_0^2(\rho)$ denotes the set of maps
 $u \in L^2(\rho)$ such that $\int u\,\dd\rho =0$.

\subsection{Linearization of the Monge-Ampère equation} Theorem~\ref{t:main} can be reformulated as stating that 
\begin{equation}
\label{eq:heuristic_goal}
\forall \mu,\nu\in \calP_2(\bbR^d), \quad \Var_\rho(\phi_\mu-\phi_\nu) \leq \frac 1 \kappa W_2^2(\mu,\nu),
\end{equation}
where  $\phi_\mu$ is the unique Brenier potential
lying in $L_0^2(\rho)$ whose gradient pushes $\rho$ forward onto~$\mu$.
To avoid any technical argument in this heuristic description of the proof, we assume
that~$\mu$ is absolutely continuous
with respect to the Lebesgue measure, 
with a smooth density, and that~$\phi_\mu$ is 
smooth and uniformly
convex.
These assumptions on $\phi_\mu$  are convenient since they imply that the Legendre transform $\phi_\mu^*$ is   smooth. We will eventually see that the smoothness and uniform convexity moduli of~$\phi_\mu$ do not appear in our bounds, so Theorem~\ref{t:main} can then be obtained through an approximation argument.
Since by definition $(\nabla\phi_\mu^*)_\# \mu=\rho$,
we deduce that $\phi_\mu^*$ is a classical solution to
the following Monge-Amp\`ere equation
$$\det(\nabla^2\phi_\mu^*) = \frac{\mu}{\rho(\nabla\phi_\mu^*)},\quad \text{on }  \bbR^d,$$
where $\rho$ and $\mu$ are identified with their respective densities.

A large body of recent work
has sought to obtain
quantitative stability bounds
for the solutions to such  
equations under perturbations of the target
measure $\mu$, motivated
by their applications to the 
optimal matching problem~\cite{caracciolo2014,ambrosio2019b,bobkov_akt,ambrosio2019,ambrosio2019c,clozeau2024annealed,goldman2024,ambrosio2024}, 
and their uses for statistical inference~\cite{manole2023,gonzalez2024, gonzalez2026empirical,gonzalez2026influence}. 
One of the outcomes of this 
literature is that, under sufficiently
strong regularity
conditions, if $\nu$ denotes
a measure close to $\mu$, also with a density,
the difference $\nabla(\phi_\mu-\phi_\nu)$
between the Brenier
maps from $\rho$ to $\mu$ and $\nu$ is
well approximated by $\nabla u$, where
$u$ is an appropriate solution to the 
partial differential equation 
\begin{equation}\label{eq:linear_pde} 
L_{\mu} u = \nu-\mu.
\end{equation}
Here,  
we define the operator 
\begin{align}\label{eq:operator_L} 
L_\mu u = -\div\big(\mu \cdot (\nabla u\circ \nabla\phi_\mu^*)),
\end{align}
acting on $u\in C_c^\infty(\R^d)$,
where we identify again absolutely continuous measures with respect to the Lebesgue measure with their densities.

One way of understanding how
$L_\mu$ relates to the Monge-Amp\`ere equation 
is through the following simple identity, which we 
state rigorously in Lemma \ref{lem:simple_identity}:
If $(\mu_t)_{t=0}^1$
is a smooth curve of probability positive densities on a uniformly convex body with smooth boundary, and $\phi_t\in L^2_0(\rho)$ is the Brenier potential from $\rho$ to $\mu_t$, then 
\begin{equation}\label{e:Lphit}
    \partial_t \mu_t = L_{\mu_t}[\partial_t\phi_t].
\end{equation}

When $\rho$ and $\mu$ are
 uniform distributions
over a set $S$ of unit volume, then  
$L_{\mu}$ is simply the negative
Laplacian on $S$, which
is the textbook
example of the linearization of 
the Monge-Amp\`ere equation~\cite{villani2003}. 
In our context where $\rho$ and $\mu$ differ,
the operator $L_\mu$
first appeared in the work
of~\cite{loeper2005}.
The general form \eqref{eq:operator_L} of $L_\mu$
resembles that of an elliptic
second-order differential operator, 
up to the fact that it contains
a composition with the map
$\nabla\phi^*$. We will 
show that the existing theory of
elliptic PDEs can nevertheless
be applied to understand the behavior
of $L_\mu$.

The identity \eqref{e:Lphit}
suggests a  natural strategy for 
deriving Brenier stability inequalities.
Given two target measures $\mu$ and $\nu$,
and an absolutely continuous
curve $(\mu_t)_{t=0}^1$ 
such that $\mu_0=\mu$ and $\mu_1=\nu$, 
one formally has
\begin{equation} \label{eq:heuristic1}
\|\phi_\mu-\phi_\nu\|_{L^2(\rho)} 
 \leq \int_0^1 \|\partial_t\phi_t\|_{L^2(\rho)}\dd t
  = \int_0^1 \|L_{\mu_t}^{-1}[\partial_t \mu_t]\|_{L^2(\rho)} \dd t,
\end{equation}
provided that 
the invertibility of $L_{\mu_t}$
can be justified. 
If one can further justify that $L_{\mu_t}$
obeys similar regularity
properties as second-order uniformly elliptic 
partial differential operators, then  
classical elliptic theory
would suggest an estimate
of the following type:
\begin{equation} 
\label{eq:heuristic2}
\|L_{\mu_t}^{-1}[\partial_t \mu_t]\|_{L^2(\rho)} 
\asymp \|\partial_t\mu_t\|_{\dot H^{-2}(\rho)}.
\end{equation}
Combining \eqref{eq:heuristic1} and \eqref{eq:heuristic2}
would formally suggest an upper bound of the form
\begin{equation}\label{eq:heuristic3}
\|\phi_\mu-\phi_\nu\|_{L^2(\rho)}
\lesssim \|\mu-\nu\|_{\dot H^{-2}(\rho)}.
\end{equation}
This suggests an
estimate which is stronger than
the one we wish to prove  in equation~\eqref{eq:heuristic_goal}. Indeed, since the Wasserstein distance
is well-known to behave infinitesimally
like a negative Sobolev norm of first order~\cite{villani2003,peyre2018},  \eqref{eq:heuristic_goal}
would formally correspond
to having $\dot H^{-1}(\rho)$ 
instead of  $\dot H^{-2}(\rho)$ on the right-hand side.

\subsection{Coercivity estimates} The heuristic
arguments leading to \eqref{eq:heuristic3}
have been made rigorous in
the special case where the 
measures satisfy
the conditions of Caffarelli's regularity
theory (see for instance Theorem 2 of~\cite{manole2023}). However, the resulting constants of 
inequality~\eqref{eq:heuristic3} depend on
the H\"older norms of the underlying measures, while we do not make the assumption that $\mu,\nu$ are absolutely continuous in the present work. Although $L_\mu$ fails to be uniformly elliptic without quantitative smoothness assumptions on $\mu$,  
we will show that equation~\eqref{eq:heuristic3}  does hold when the $\dot H^{-2}(\rho)$ norm
is replaced by the $\dot H^{-1}(\rho)$ norm, 
even without such 
assumptions on $\mu$ and $\nu$.

Our main tool for establishing this fact is the following 
estimate for $L_\mu$: 
\begin{equation}\label{eq:Akappa_heuristic}
 \forall u\in \calC_c^\infty(\bbR^d)\cap L_0^2(\rho), 
\quad \|u\|_{L^2(\rho)}^2 \leq \frac 1 \kappa \| L_\mu u\|_{\dot H^{-1}(\mu)}^2.
\end{equation}
Equation~\eqref{eq:Akappa_heuristic} is reminiscent of   $H^1$ regularity estimates
for elliptic operators; we will refer to it
as a {\it coercivity} estimate. 
The key feature of~\eqref{eq:Akappa_heuristic} is that the constant $\kappa$ does not depend
on the regularity of $\mu$. Similar estimates have
been proven recently, but only 
under the stronger assumption that both $\rho$ and $\mu$ are  uniformly log-concave~\cite{han2026sharp}.

To establish coercivity estimates like~\eqref{eq:Akappa_heuristic}, we rely on generalizations of Bochner's formula. 
To gain some intuition for its relevance here, let us 
formally consider the regime  where $\mu = \rho$, 
in which case equation~\eqref{eq:Akappa_heuristic}
boils down 
to the Poincaré inequality 
$$
\forall u\in \calC_c^\infty(\bbR^d)\cap L_0^2(\rho), \quad \|u\|_{L^2(\rho)}\leq \kappa^{-1/2}\|\Delta_\rho u\|_{\dot{H}^{-1}(\rho)}=\kappa^{-1/2}\|\nabla u\|_{L^2(\rho)}.
$$
A classical proof of this inequality goes through Bochner's integrated formula
\begin{equation}\label{e:bochner}
\int(\Delta_\rho u)^2\mathrm{d}\rho = \int\|\nabla^2 u\|_{F}^2\, \mathrm{d}\rho +\int\langle\nabla^2V\nabla u,\nabla f\rangle \,\mathrm{d}\rho.
\end{equation}
Let us recall this proof. Denote by $\lambda>0$ the first non-trivial eigenvalue of $-\Delta_\rho$, and let $h$ be an associated eigenfunction. 
Using \eqref{e:bochner} and the fact that $\|\nabla^2 h(x)\|_{F}^2\geq 0$, we get 
$$
\lambda^2\|h\|_{L^2(\rho)}^2\geq \kappa \|\nabla h\|_{L^2(\rho)}^2=\kappa\int (-\Delta_\rho h)h\, \dd\rho=\kappa\lambda \|h\|_{L^2(\rho)}^2.
$$
Therefore $\lambda\geq \kappa$, which implies that for any $f\in W^{1,2}(\rho)\cap L^2_0(\rho)$,
$$
\|\nabla f\|_{L^2(\rho)}^2 =\int (-\Delta_\rho f)f\dd\rho \geq  \lambda\|f\|_{L^2(\rho)}^2\geq \kappa\|f\|_{L^2(\rho)}^2.
$$

In the general case where $\mu\neq\rho$, the operator $L_\mu$ is not simply a weighted Laplacian: it contains a composition with $\nabla\phi^*$ which makes it in particular non-symmetric. To handle this
more general case, we first prove a further reduction, showing
that the coercivity
estimate~\eqref{eq:Akappa_heuristic} is implied by the coercivity of a distinct
operator $E_\mu$, 
which can be thought
of as the pullback of $L_\mu$ 
through the map $\nabla \phi_\mu^*$. This operator already appeared  in the study of Poincaré inequalities for log-concave measures, notably in works by Kolesnikov, Klartag and Milman (see Remark \ref{r:Emu} for references). 
In contrast to $L_\mu$, 
this pullback
operator is symmetric,
and our proof of its coercivity
relies on a generalization  of Bochner's formula specifically suited to composition with gradients (see Lemma \ref{l:genbochner} below).

In order to make these various
arguments precise, we 
first truncate $\rho$ to a bounded set. Doing
so allows us to appeal to
the 
regularity theory of optimal transport
maps due to Caffarelli~\cite{caffarelli1992a,caffarelli1992,caffarelli1996} and Urbas~\cite{urbas1988regularity,urbas1997second}, 
and the resulting H\"older
regularity in time of curves
of Brenier potentials~\cite{gonzalez2024}. 
For coercivity estimates
and generalizations of Bochner's identity, we then need to work with appropriate boundary conditions. 

We turn in the next sections to 
rigorously proving our main results.
We begin by   establishing
the reduction from coercivity estimates to 
Brenier stability inequalities, in a sufficiently
general setting so as to be applicable
to proving both Theorems~\ref{t:main} and~\ref{t:convexset}.

\section{Reduction to coercivity estimates}\label{s:reductocoercest}
In this section, we establish Brenier stability
inequalities under 
appropriate coercivity conditions
on the operator $L_\mu$.  
Our results are stated under 
strong regularity conditions, which will later be relaxed
by means of regularization and truncation 
arguments.

To fix notation for the remainder
of this section, let $\calQ$ denote the set
of   uniformly convex bodies
in $\bbR^d$ with $\calC^\infty$ boundary. 
Given $M \in \calQ$, we denote
by $n_M$ the outward-pointing unit
normal vector
of $M$. 
Given $\gamma >0$,
  let  $\calF(M;\gamma)$ denote set of
  Borel probability measures $\sigma$  on $M$
  which are absolutely continuous with respect to the Lebesgue measure
on $M$,
with densities lying in $\calC^\infty(M)$
and satisfying
$$\|\sigma\|_{\calC^2(M)} \leq \gamma,\quad 
\text{and}\quad \sigma \geq \gamma^{-1} ~~~\text{over } M.$$
Furthermore, let $\calF(Y) = \bigcup_{\gamma > 0} \calF(Y;\gamma)$.
For any   measure $\sigma \in \calF(Y)$, we further define
the homogeneous Sobolev norm
\begin{equation} \label{eq:sobolev_norms}
\|f\|_{\dot H^1(\sigma)} = \|\nabla f \|_{L^2(\sigma)}
\end{equation}
for all $f \in \calC^\infty(M)\cap L_0^2(\sigma)$, 
as well as its dual norm
\begin{equation}\label{e:H-1def}
\|h\|_{\dot H^{-1}(\sigma)} = 
\sup\left\{ \int_M h(x)g(x)\dd x: g \in \calC^\infty(M)\cap L_0^2(\sigma),
\|g\|_{\dot H^1(\sigma)}\leq 1 
\right\},
\end{equation}
for all $h\in $
Furthermore, we denote by $\dot H^1(\sigma)$ and $\dot H^{-1}(\sigma)$
the completion of $\calC^\infty(M)\cap L_0^2(\sigma)$ 
with respect to the 
norms in~\eqref{eq:sobolev_norms} and~\eqref{e:H-1def}. 
Our assumptions on $\sigma$  
imply that 
$\dot H^1(\sigma) = W^{1,2}(\sigma)\cap L_0^2(\sigma)$.

\subsection{The linearization operator $L_\mu$}\label{s:Lmuop}

Consider
two   convex bodies
 $K,Y \in \calQ$,
 and  given $\gamma > 0$,   fix a measure $\rho \in \calF(K;\gamma)$.
Under these conditions, the 
classical 
boundary regularity
theory of optimal transport~\cite{caffarelli1996,urbas1997second} implies that
there exists a constant $\lambda = \lambda(\gamma,K,Y) > 0$ such that
for any 
$\mu \in \calF(Y;\gamma)$, 
the unique mean-zero Brenier potential $\phi_\mu$ 
whose gradient pushes $\rho$ forward onto $\mu$   satisfies
\begin{equation} 
\label{eq:Cinfty_potentials} 
\phi_\mu \in \calC^\infty(K), \phi_\mu^* \in \calC^\infty(Y),
\end{equation} 
so that 
the map $\nabla\phi_\mu:K\to Y$
defines a $\calC^\infty$-diffeomorphism,
and furthermore, 
\begin{equation}
\label{e:hessphi}
\frac 1 \lambda I_d 
\preceq 
\nabla^2\phi_\mu(x),\nabla^2\phi_\mu^*(y)\preceq 
\lambda I_d, 
\quad x \in {\rm int}(K), y\in {\rm int}(Y).
\end{equation}
Under these conditions,  the operator
$$L_\mu u (y) = -\div\big(\mu(y)\nabla u (\nabla\phi_\mu^*(y))\big),\quad y \in {\rm int}(Y)$$
  is well-defined
in the classical sense for all $u \in \calC^2(K).$
Our main motivation
for defining $L_\mu$ is the following
identity, announced in equation~\eqref{e:Lphit}, which we now
prove rigorously.
\begin{lemma}\label{lem:simple_identity}
Let $\rho \in \calF(K)$,
and for some $\gamma > 0$, 
let $\mu_t\in \calF(Y;\gamma )$
for all $t\in [0,1]$.
Suppose further that 
$t \in [0, 1] \mapsto \mu_t \in \calC^{2}(Y)$ is a $\calC^1$ map between Banach spaces. Let $\phi_t\in L^2_0(\rho)$ denote the Brenier potential from $\rho$ to $\mu_t$. Then, 
$$\partial_t \mu_t = L_{\mu_t}[\partial_t\phi_t].$$
\end{lemma}
Variants
of this result 
already exist in the literature (e.g.~\cite{gonzalez2026empirical}),
but we were unable to find
a reference with this precise formulation, 
so we provide a  proof~below.
\begin{proof}[Proof of Lemma~\ref{lem:simple_identity}]
Notice that $L_{\mu_t}[\partial_t\phi_t]$
is well-defined: Indeed, thanks to the uniform convexity of $K$,
Theorem 3.2 of~\cite{gonzalez2024}
implies that the mapping $t\in [0,1]\mapsto \phi_t\in \calC^{2}(K)$ is $\calC^1$ as a curve in $\calC^{2}(K)$.
Now, for any test function 
$f \in \calC_c^\infty(\bbR^d)$, 
the pushforward condition $\mu_t = (\nabla\phi_t)_\# \rho$ implies that
$$\int_Y f \dd\mu_t = \int_K f(\nabla\phi_t)\dd\rho$$
By 
taking a derivative over $t$, we obtain  
\begin{align}
\nonumber 
    \frac \dd {\dd t} \int_Y f\dd\mu_t 
 &= \int_K \langle \nabla f(\nabla\phi_t), \partial_t\nabla\phi_t\rangle \dd\rho\\
 &=  - \int_Y  f\cdot  \div\big(\mu_t \cdot (\partial_t\nabla\phi_t)\circ (\nabla\phi_t^*)
\big) + \int_{\partial Y} f \mu_t \big\langle (\partial_t\nabla \phi_t) \circ \nabla\phi_t^*, n_Y\big\rangle.
\label{eq:step_Lmu}
 \end{align}
 Let us show that
\begin{equation}\label{e:potatbdry}
\langle (\partial_t \nabla\phi_t)\circ \nabla\phi_t^*,n_Y\rangle = 0
\end{equation}
 over $\partial Y$, which will conclude the proof of the claim. This boundary
 condition has been previously 
 identified
 for instance by~\cite{gonzalez2024}. 
To prove it, 
recall that $\nabla\phi_t$
is a diffeomorphism 
between $K$ and $Y$, 
thus one must have $\nabla\phi_t(\partial K) = \partial Y$ for all $t \in [0,1]$. 
Therefore, given any defining function $h \in \calC^\infty(\bbR^d)$ 
for the set $Y\in\calQ$, satisfying
$${\rm int}(Y) =\{h < 0\}, 
\quad \partial Y = \{h=0\},\quad 
\nabla h \neq 0~~\text{on } \partial Y,$$
we must have
$$h(\nabla\phi_t(x)) = 0,\quad x \in \partial K.$$
Differentiating implies
$$\langle \nabla h(\nabla\phi_t), \partial_t \nabla\phi_t\rangle = 0,\quad 
\text{on } \partial K. $$
Since $\partial Y=\{h=0\}$, the restriction of $\nabla h$ to $\partial Y$
is parallel to the outward-pointing
unit vector $n_Y$. We deduce
$$\langle n_Y(\nabla\phi_t),\partial_t\nabla\phi_t\rangle = 0,\quad \text{on } \partial K,$$
which is equivalent to 
$$\langle n_Y,(\partial_t\nabla\phi_t)\circ (\nabla\phi_t^*)\rangle = 0,\quad \text{on } \partial Y.$$
This proves \eqref{e:potatbdry},  
and the claim follows.
\end{proof}

\subsection{A general criterion for Brenier stability
inequalities}
Given a convex body $K\in \calQ$ and a measure
$\rho \in \calF(K)$, we
say that $\rho$ satisfies the {\it coercivity condition $A(\kappa)$} 
for some constant $\kappa > 0$ if the following
assertion holds.
\vspace{0.05in}
\begin{enumerate}
\item[$A(\kappa)$] 
For all $Y\in \calQ$, $\mu\in\calF(Y)$,
and $u \in \calC^2(K)\cap L_0^2(\rho)$
satisfying the boundary
condition
\begin{equation}
\forall y\in \partial Y, \quad \langle \nabla u(\nabla\phi_\mu^*(y)),
n_Y(y)\rangle = 0,
\end{equation}
there holds
\begin{equation}\label{e:coerciveLmu}
\|u\|_{L^2(\rho)} \leq \kappa^{-1/2} \|L_\mu u\|_{\dot{H}^{-1}(\mu)}.
\end{equation}
\end{enumerate} 
\vspace{0.05in}
Under the  general coercivity condition $A(\kappa)$, we establish the following
 stability bound.
\begin{theorem}\label{t:general_reduction}
Let $K \in \calQ$ be a convex body, and $\rho \in \calF(K)$.
Assume   that~$\rho$   satisfies condition $A(\kappa)$ for some $\kappa > 0$. 
Then, $\rho$ satisfies a Brenier
stability inequality with constant 
$$C_W(\rho) \leq \frac 1 \kappa.$$
\end{theorem}

Theorem~\ref{t:general_reduction} shows that the coercivity
condition $A(\kappa)$
is  sufficient   for proving
a Brenier stability inequality
for a compactly-supported  measure $\rho$. 
Importantly, we emphasize that
the upper
bound on the constant $C_W(\rho)$
does not depend on 
any regularity properties of $\rho \in \calF(K)$.

The remainder of this subsection is devoted
to the proof of Theorem~\ref{t:general_reduction}. 
We begin with the following
Lemma, which 
builds upon Lemma~\ref{lem:simple_identity}
to obtain an upper
bound on the variance of differences
of Brenier potentials.
 
\begin{lemma}
\label{lem:main_smooth}
Let $\rho$ be a measure
satisfying the same conditions
as in the statement of Theorem~\ref{t:general_reduction},
and let $\mu_t \in \calF(Y;\gamma)$
denote a curve of measures
indexed by $t \in [0,1]$
for some $\gamma > 0$. 
Suppose that the curve
$t \in [0, 1] \mapsto \mu_t \in \calC^{2}(Y)$ is a $\calC^1$ mapping between Banach spaces.
Suppose further that there exists
a velocity field $(v_t)_{t=0}^1$ that classically satisfies
the  
continuity equation with no-flux boundary condition,  
\begin{align*} 
\partial_t \mu_t + \div(v_t \mu_t) &= 0,\quad \text{in } Y,\\
\langle v_t,n_Y\rangle &= 0,\quad \text{on } \partial Y.
\end{align*}
Furthermore, let $\phi_0,\phi_1\in L^2_0(\rho)$ denote Brenier potentials from $\rho$ to $\mu_0$ and $\mu_1$ respectively. Then
\begin{equation}\label{e:main_smooth}
\Var_\rho(\phi_1-\phi_0)\leq \frac1\kappa \int_0^1\|v_t\|^2_{L^2(\mu_t)} \dd t. 
\end{equation}
\end{lemma}
\begin{proof}[Proof of Lemma \ref{lem:main_smooth}]
Denote by $\phi_t\in L^2_0(\rho)$ the Brenier potential such that $(\nabla\phi_t)_\#\rho=\mu_t$.  
Similarly to the proof of Lemma \ref{lem:simple_identity},
thanks  to the uniform convexity of $\partial K$,
Theorem 3.2 of~\cite{gonzalez2024}
implies that the mapping $t\mapsto \phi_t\in \calC^{2}(K)$ is $\calC^1$ as a curve in $\calC^{2}(K)$.
Furthermore, by the same
argument as in equation~\eqref{e:potatbdry},
the time
derivative $u:=\partial_t\phi_t\in \calC^2(K)\cap L_0^2(\rho)$ satisfies
the boundary condition
$$\langle \nabla u \circ \nabla\phi_t^*, n_Y\rangle = 0\quad \text{on } \partial Y,$$
thus all the assumptions
are satisfied to apply
the coercivity condition~$A(\kappa)$  with this choice of test function $u$, and 
with 
$\phi=\phi_t$, to obtain 
\begin{align*} 
\|\partial_t \phi_t\|_{L^2(\rho)} 
 &\leq \kappa^{-1/2} \|L_{\mu_t}[\partial_t\phi_t]\|_{\dot H^{-1}(\mu_t)}.
 \end{align*}
 Furthermore,  $(\mu_t)_{t=0}^1$ satisfies the assumptions of Lemma \ref{lem:simple_identity}, thus we further obtain
\begin{align}
\label{eq:divstep}
\|\partial_t \phi_t\|_{L^2(\rho)} 
 \leq \kappa^{-1/2}  \|\partial_t \mu_t\|_{\dot H^{-1}(\mu_t)}
 =\kappa^{-1/2}  \|\div (v_t\mu_t)\|_{\dot H^{-1}(\mu_t)}.
\end{align}  
Let us observe that
\begin{equation}\label{e:transportH-1L2}\|\div(v_t\mu_t)\|_{\dot H^{-1}(\mu_t)}
\leq \|v_t\|_{L^2(\mu_t)}.\end{equation}
Indeed, take any $g\in C^\infty(M)\cap L^2_0(\mu_t)$ with $\|\nabla\varphi\|_{L^2(\mu_t)} \leq 1$. By Cauchy–Schwarz,
$$\langle \div(v_t\mu_t), g\rangle = -\int \nabla g\cdot v_t\,\dd\mu_t \leq \|\nabla g\|_{L^2(\mu_t)}\, \|v\|_{L^2(\mu_t)} \leq \|v_t\|_{L^2(\mu)}.$$
Taking the supremum over $g$ gives \eqref{e:transportH-1L2}.
Combining this bound
with equation~\eqref{eq:divstep}, 
we thus obtain
\begin{align*} 
\|\phi_1-\phi_0\|_{L^2(\rho)}
\leq \int_0^1 \left\|  \partial_t \phi_t\right\|_{L^2(\rho)}\dd t
\leq \kappa^{-1/2} \left(\int_0^1\|v_t\|^2_{L^2(\mu_t)} \dd t\right)^{1/2}
\end{align*} 
which concludes the proof.
\end{proof}

We finally conclude the proof of Theorem \ref{t:general_reduction} using the fact that Lemma \ref{lem:main_smooth}, although proved under strong regularity assumptions on the curve of measures $(\mu_t)_{t=0}^1$, has a bound which does not depend on these assumptions. 
\begin{proof}[Proof of Theorem \ref{t:general_reduction}]
Let $\mu_i=(\nabla\phi_i)_\#\rho$ and assume for a moment that $\mu_0,\mu_1$ are compactly supported. Let $Y\subset\R^d$ be a closed ball 
centered at the origin, and containing in its interior the supports of both $\mu_0$ and $\mu_1$. Consider $\pi$ an optimal coupling between $\mu_0$ and $\mu_1$, and denote by $\nu_t$ the law of the random variable $(1-t)y+tz$ when $(y,z)\sim \pi$. 
Let $\chi$ be a smooth non-negative function supported in the unit ball, with $\int \chi \,\dd x=1$. Set $\chi_\delta=\delta^{-d}\chi(\delta^{-1}\cdot)$. For $\delta>0$ small enough, the regularized density
$$
\mu_{t}^{(\delta)}(x)=\chi_\delta\star \nu_t(x)=\int \chi_\delta(x-(1-t)y-tz) \,\dd\pi(y,z)
$$
is supported in the interior of $B$ for any $t\in[0,1]$. Moreover, it satisfies 
$$
\partial_t\mu_t^{(\delta)}+\div(j_t^{(\delta)})=0
$$
for 
$$
j_t^{(\delta)}(x)=\int (z-y)\chi_\delta(x-(1-t)y-tz)\,\dd\pi(y,z).
$$
We define the velocity by $v_t^{(\delta)}=\frac{j_t^{(\delta)}}{\mu_t^{(\delta)}}$ 
on $\{\mu_t^{(\delta)}>0\}$, and $v_t^{(\delta)}=0$ elsewhere.
Then, denote by $\omega$ the uniform probability measure on $Y$. For $0<\varepsilon<1$, let
$$
\mu_t^{(\delta,\varepsilon)}=(1-\varepsilon)\mu_t^{(\delta)}+\varepsilon \omega
$$
whose density is bounded below on $Y$ uniformly over $t\in[0,1]$. In particular, $\mu_t^{(\delta,\varepsilon)}\in \mathcal{F}(Y;\gamma)$ for some $\gamma>0$. Moreover, $t\mapsto \mu_t$ is a $\calC^1$ map from $[0,1]$ to $\calC^{2}(Y)$, and it satisfies
$$
\partial_t\mu_t^{(\delta,\varepsilon)}+\div(j_t^{(\delta,\varepsilon)})=0
$$
for $j_t^{(\delta,\varepsilon)}=(1-\varepsilon)j_t^{(\delta)}$. We 
then define the velocity field~$v_t^{(\delta,\varepsilon)}=\frac{j_t^{(\delta,\varepsilon)}}{\mu_t^{(\delta,\varepsilon)}}$.
Notice that the support of~$j_t^{(\delta)}$ 
is compactly contained in the interior of $Y$ for all~$\delta$
sufficiently small, thus the velocity
field~$v_t^{(\delta,\varepsilon)}$ vanishes on $\partial Y$. 
In particular, it satisfies  
\begin{equation} 
\langle v_t^{(\delta,\varepsilon)},n_Y\rangle = 0,
\quad \text{on }\partial Y.
\end{equation}
Denoting by $\phi_i^{(\delta,\varepsilon)}\in L^2_0(\rho)$ the Brenier potential from $\rho$ to $\mu_i^{(\delta,\varepsilon)}$ for $i=0,1$,
the preceding two displays allow us to apply 
Lemma \ref{lem:main_smooth}, to obtain
\begin{equation}\label{e:variancedeltaeps}
\Var_\rho(\phi_1^{(\delta,\varepsilon)}-\phi_0^{(\delta,\varepsilon)})\leq \frac1\kappa \int_0^1 \|v_t^{(\delta,\varepsilon)}\|_{L^2(\mu_t^{(\delta,\varepsilon)})}^2 \,\dd t.
\end{equation}
We then take the limit $\delta,\varepsilon\rightarrow 0$. For this, we first recall that the map from $(\mathcal{P}_2(\R^d),W_2)$ to $L^2(\rho,\R^d)$ which to $\mu$ associates the Brenier map from $\rho$ to $\mu$ is continuous. Thus, by the Poincaré inequality for $\rho$, which holds
since $\rho\in\calF(K)$, the map from $(\mathcal{P}_2(\R^d),W_2)$ to $L^2_0(\rho)$ which to $\mu$ associates the Brenier potential from $\rho$ to $\mu$ is continuous. Since $W_2(\mu_i^{\delta,\varepsilon},\mu_i)\rightarrow 0$ for $i\in\{0,1\}$, we get that the left-hand side of \eqref{e:variancedeltaeps} tends to $\Var_\rho(\phi_1-\phi_0)$. Regarding the right-hand side let us prove that
\begin{equation}\label{e:vtdeltaeps}
\int_0^1 \|v_t^{(\delta,\varepsilon)}\|_{L^2(\mu_t^{(\delta,\varepsilon)})}^2 \,\dd t \underset{\delta,\varepsilon\rightarrow 0}{\longrightarrow} \int |z-y|^2\, \dd\pi(y,z)=W_2^2(\mu_0,\mu_1)
\end{equation}
since once this is proven, \eqref{e:varkappaW2} follows from letting $\delta,\varepsilon\rightarrow 0$ in \eqref{e:variancedeltaeps}.
To prove \eqref{e:vtdeltaeps}, we first make the elementary observation that
$$
\int_0^1 \|v_t^{(\delta,\varepsilon)}\|_{L^2(\mu_t^{(\delta,\varepsilon)})}^2 \,\dd t  \rightarrow \int_0^1 \|v_t^{(\delta)}\|_{L^2(\mu_t^{(\delta)})}^2 \,\dd t $$
as $\varepsilon\rightarrow 0$. Then, the Benamou-Brenier formula \cite{benamou2000computational} yields 
$$
\int_0^1\|v_t^{(\delta)}\|^2_{L^2(\mu_t^{(\delta)})}\,\dd t\geq  W_2^2\big(\mu_0^{(\delta)},\mu_1^{(\delta)}\big)
$$
and the right-hand side tends to $W_2^2(\mu_0,\mu_1)$ as $\delta\rightarrow 0$. For the converse inequality, we notice that by the Cauchy-Schwarz inequality,
$$\|v_t^{(\delta)}\|^2_{L^2(\mu_t^{(\delta)})}\le \iint |z-y|^2\chi_\delta(x-(1-t)y-tz)\,\dd x\,\dd\pi(y,z)=\int|z-y|^2\,\dd\pi = W_2^2(\mu_0,\mu_1)$$
and integrating over $t\in[0,1]$ we finally conclude the proof of \eqref{e:vtdeltaeps}.

Let us finally address the case where $\mu_0,\mu_1$ are not necessarily compactly supported. For $i\in\{0,1\}$, consider $(\mu_i^{(n)})_{n\in\N}$ a sequence of probability measures which converges in $W_2$ norm to $\mu_i$, and such that each $\mu_i^{(n)}$ is compactly supported. Denote by $\phi_i^{(n)}$ the Brenier potential from $\rho$ to $\mu_i^{(n)}$. We already know by the previous step that 
$$
\Var_\rho(\phi_1^{(n)}-\phi_0^{(n)})\leq \frac1\kappa W_2^2(\mu_0^{(n)},\mu_1^{(n)}).
$$
The right-hand side converges to $\frac1\kappa W_2^2(\mu_0,\mu_1)$ by the triangle inequality. The left-hand side converges to $\Var_\rho(\phi_1-\phi_0)$ once again by the continuity of the map  from $(\mathcal{P}_2(\R^d),W_2)$ to $L^2_0(\rho)$ which to $\mu$ associates the Brenier potential from $\rho$ to $\mu$. Hence, \eqref{e:varkappaW2} holds, which concludes the proof.
\end{proof}

\subsection{Reduction to a coercivity estimate on an auxiliary operator $E_\mu$} \label{s:reduccoerciveEphi}
In this section, we show that the coercivity estimate on $L_\mu$ appearing as equation~\eqref{e:coerciveLmu} in condition $A(\kappa)$ can be deduced from a coercivity estimate on an auxiliary operator $E_\mu$ which has the advantage of being symmetric, and which is obtained by pulling back $L_\mu$ to the source domain $K$. 

We again fix convex bodies $K,Y\in \calQ$, 
and measures $\rho \in \calF(K;\gamma)$ and 
$\mu \in \calF(Y;\gamma)$ for
some $\gamma > 0$. The operator $L_\mu$ in \eqref{eq:operator_L} 
is generally not self-adjoint due
to the fact that it contains a composition
with the map $\nabla\phi^*$. Recalling that $\nabla\phi^*\circ\nabla \phi={\rm id}$ (which is a direct consequence of~\eqref{e:convexconj}), it is therefore natural to consider the operator $u\mapsto L_\mu u \circ (\nabla\phi)$. As observed in Lemma~8 of \cite{manole2023} (see also \eqref{e:EphiandLphi} below), this manipulation naturally leads to the operator
\begin{equation}\label{eq:operator_E} 
E_\mu u = - \div_\rho((\nabla^2\phi)^{-1}\nabla u),
\end{equation}
where we recall that 
$\div_\rho(u) = \div(\rho u)/\rho$ denotes the weighted divergence with respect to $\rho$.

\begin{remark}[The Operator $E_\mu$]\label{r:Emu}
We have
seen  that the operator $L_\mu^{-1}[\nu-\mu]$ is formally the linearization of the difference $\phi_\nu-\phi_\mu$ between sufficiently-regular Brenier potentials.
Likewise, 
the operator $E_\mu$ 
arises in the linearization 
of their convex conjugates, $\phi_\mu^*-\phi_\nu^*$, and 
therefore plays an important
role in the 
aforementioned literature
on linearization of
Monge-Amp\`ere equations. 
On the other hand, 
$E_\mu$ also arises
in a separate
line of work in convex
geometry. When transporting a log-concave density $\rho$ to another one $\mu$ via a transport map $\nabla \phi$, it is natural to study the ``weighted Riemannian manifold" $(\R^d, \nabla^2\phi, \rho)$, meaning that $\R^d$ is endowed with the Hessian metric $\nabla^2\phi$ and with the log-concave probability measure $\rho$. The natural ``Laplace-type" operator on this manifold is $E_\mu$. This perspective has been developed notably in \cite{kolesnikov2014,klartag2014logarithmically,klartag2015eigenvalue,kolesnikov2016riemannian}, and most recently in \cite{han2026sharp}. Most of this literature focuses on the case where both $\rho$ and $\mu$ are log-concave, while we are dealing  merely with a smooth density $\mu$ on a compact convex domain. While the coercivity estimate on $E_\mu$ given by the condition $B(\kappa)$ below was not known in this general set-up, it is worth mentioning \cite[Proposition 2.4]{han2026sharp} where the condition $A(\kappa)$ has been established in the framework where both $\rho$ and $\mu$ are strongly log-concave.
\end{remark}
Unlike $L_\mu$, the operator $E_\mu$   is symmetric with respect to the $L^2(\rho)$ scalar product,
and uniformly elliptic under the conditions
we have placed on $K,Y,\rho,\mu$,
in view of equation~\eqref{e:hessphi}.  
We can thus easily prove the 
existence of solutions to equations of the form
$E_\mu u = f$. Consider the boundary value problem
\begin{equation}
\label{eq:bvp}
    \begin{aligned}
        E_\mu u&= f, \quad \text{in } {\rm int}(K), \\
        \langle (\nabla^2\phi)^{-1}\nabla u, n_K\rangle &= 0, \quad \text{on $\partial K$}
    \end{aligned}
\end{equation}
where we recall that $K$ is the support of
$\rho$ and $n_K$ is the outward-pointing
unit normal vector to $K$. 
In order to analyze this problem, define
the following bilinear form for $u,v \in \dot H^1(\rho)$
\begin{equation}\label{e:phiscal}
\langle u,v\rangle_\phi =  
\int \langle (\nabla^2\phi)^{-1}\nabla u, \nabla v\rangle\, \dd\rho,
\end{equation}
Notice that when $u,v \in \calC^\infty(\bbR^d)$ and $u$ satisfies equation~\eqref{eq:bvp},
we   have $\langle u,v\rangle_\phi = \int (E_\mu u)  v \, \dd\rho$. 
Therefore, for any $f \in L_0^2(\rho)$, 
we say that $u \in \dot H^1(\rho)$ is
a (weak) solution to the Neumann problem~\eqref{eq:bvp} 
provided  
\begin{equation}\label{e:phiscal2}
\forall v \in \calC^\infty(\bbR^d), \qquad \langle u,v\rangle_\phi = \langle f,v\rangle_{L^2(\rho)}.\end{equation}
By density, this identity then holds for any $v\in H^1(\rho)$.
\begin{lemma}\label{lem:weak_solutions}
Let $K,Y \in \calQ$ be convex bodies.
Furthermore, let $\rho\in\calF(K)$
and $\mu \in \calF(Y)$.
Then, for any $f \in L_0^2(\rho)$, there exists
a unique weak solution $u \in \dot H^1(\rho)$ of
the Neumann problem~\eqref{eq:bvp}. Moreover, the application $f\mapsto u$ is continuous from $L^2_0(\rho)$ to $\dot H^1(\rho)$.
\end{lemma}
\begin{proof} 
Under the assumptions
we have placed on $\rho$ and $\mu$, 
the map $\phi$ satisfies equation~\eqref{e:hessphi},
thus for any $u \in \dot H^1(\rho)$
$$\langle u,u\rangle_\phi \geq  \lambda^{-1}\|\nabla u\|^2_{L^2(\rho)} = \lambda^{-1}\|u\|^2_{\dot  H^1(\rho)}$$
where the second inequality follows from the Poincaré inequality (assumed in $A(\kappa)$).
Also by equation~\eqref{e:hessphi}, we have
$$|\langle u,v\rangle_\phi| \leq \lambda \|u\|_{\dot H^1(\rho)}\|v\|_{\dot H^1(\rho)},
\quad \text{for all } u,v \in \dot H^1(\rho).$$
The claim then follows from the Lax-Milgram theorem.
\end{proof}

The next lemma asserts than an appropriate coercivity estimate on $E_\mu$ implies the sought coercivity estimate on $L_\mu$. We will say that the measure $\rho$ 
satisfies the {\it coercivity condition $B(\kappa)$} 
for some constant $\kappa > 0$ if the following
assertion holds.
\begin{enumerate}
\item[$B(\kappa)$] 
For all $Y\in \calQ$, $\mu \in \calF(Y)$,
and $f\in L^2_0(\rho)$, there holds 
\begin{equation}\label{e:loweronEphi}
\|(\nabla^2\phi_\mu)^{-1}\nabla u\|_{L^2(\rho)}\leq \kappa^{-1/2}\|f\|_{L^2(\rho)}
\end{equation}
where $u\in \dot H^1(\rho)$ is the unique weak solution of $E_\mu u=f$.
\end{enumerate}
\begin{lemma}\label{lem:Ephi}
The coercivity condition $B(\kappa)$ implies the coercivity condition $A(\kappa)$.
\end{lemma}
\begin{proof}[Proof of Lemma \ref{lem:Ephi}]
Let $Y \in \calQ$ and $\mu \in \calF(Y)$ be
arbitrary. 
Let $f \in \calC^2(K)\cap L_0^2(\rho)$
satisfy the boundary condition
\begin{align}\label{eq:bdrycond3} 
\langle \nabla f \circ \nabla\phi_\mu^*, n_Y\rangle = 0,
\quad \text{on } \partial Y,
\end{align} 
and let $u \in \dot H^1(\rho)$ be the unique weak solution
to the boundary value problem~\eqref{eq:bvp},
which exists thanks to our assumption
on $\phi_\mu$ and Lemma~\ref{lem:weak_solutions}.  
Since $\nabla\phi_\mu^* \circ \nabla\phi_\mu={\rm id}$, we have 
\begin{equation}\label{e:hessphihessphi*}
\nabla^2\phi_\mu^* \circ \nabla\phi_\mu = (\nabla^2\phi_\mu)^{-1}.
\end{equation} 
According to \eqref{e:phiscal}-\eqref{e:phiscal2}, and using \eqref{e:hessphihessphi*}, we get 
\begin{equation} \label{e:EphiandLphi}
\begin{aligned}
\|f\|_{L^2(\rho)}^2
 &= \int \langle \nabla f, (\nabla^2\phi_\mu)^{-1} \nabla u\rangle \, \dd\rho \\
 &= \int \langle \nabla f\circ \nabla\phi_\mu^*, ((\nabla^2\phi_\mu)^{-1}\circ \nabla\phi_\mu^*) (\nabla u \circ \nabla\phi_\mu^*)\rangle \, \dd\mu \\
 &= \int \langle \nabla f\circ \nabla\phi_\mu^*,   \nabla (u \circ\nabla\phi_\mu^*)\rangle \, \dd\mu \\
 &= \langle L_\mu f, u\circ \nabla\phi_\mu^*\rangle_{L^2(\bbR^d)},
\end{aligned}
\end{equation}
where the integration by parts
in the final equality
is justified due to equation~\eqref{eq:bdrycond3}.
Since $\int L_\mu f \, \dd x=0$, it follows from equation~\eqref{e:EphiandLphi} and 
\eqref{eq:sobolev_norms} that
\begin{equation}\label{e:fg}
\|f\|_{L^2(\rho)}^2 \leq \|L_\mu f\|_{\dot H^{-1}(\mu)} 
       \|u\circ \nabla\phi_\mu^*\|_{\dot H^1(\mu)}
       = \|L_\mu f\|_{\dot H^{-1}(\mu)} 
       \|g\|_{L^2(\rho)}
\end{equation}
where $g = (\nabla^2\phi_\mu)^{-1}\nabla u$. The last equality in \eqref{e:fg} is a consequence of the identity $\nabla(u\circ\nabla\phi_\mu^*)=(\nabla^2 \phi_\mu^*)(\nabla u \circ \nabla \phi_\mu^*)$ and \eqref{e:hessphihessphi*}.
Lemma \ref{lem:Ephi} provides the bound
\begin{equation} 
\label{eq:pf_subelliptic_suff}
\|g\|_{L^2(\rho)}\leq \kappa^{-1/2}\|f\|_{L^2(\rho)}.
\end{equation}
Plugging into \eqref{e:fg}, this concludes the proof of Lemma \ref{lem:Ephi}.
\end{proof}

The goal of the next two sections is to prove
the coercivity condition $B(\kappa)$ 
under appropriate assumptions on the source measure
$\rho$, and to deduce our main results.  Section \ref{s:coercive_ulc} is devoted to proving $B(\kappa)$ when $\rho$ is uniformly log-concave on a convex body, and proving Theorem \ref{t:main} by truncation arguments. Then, Section \ref{s:coercive_cvx} 
establishes condition $B(\kappa)$ when $\rho$ is any log-concave on a convex body, but with a worse constant, and Theorem \ref{t:convexset} follows.

\section{Proof of Theorem \ref{t:main} and Corollary \ref{c:bures}}
\label{s:coercive_ulc} 

We assume in the first part of this section that $K$ is a uniformly convex body with $\calC^\infty$ boundary, 
which we recall is abbreviated by $K\in \calQ$. 
Furthermore, we assume that the 
measure $\rho$ satisfies
\begin{equation} \label{e:unif_log_concave}
\dd\rho(x)=e^{-V(x)}\dd x, \quad \text{with }
\nabla^2 V(x)\succeq \kappa I_d
~\text{ for any } x\in K
\end{equation} 
for some $\kappa>0$.
Our first goal is to prove the following result, which together with regularization arguments implies Theorem \ref{t:main} (see the end of this section).
\begin{theorem}\label{t:uniflogconcBkappa}
If $K\in \calQ$, and $\rho$ satisfies \eqref{e:unif_log_concave}, then condition $B(\kappa)$ holds. 
\end{theorem}

In the sequel we work under the assumptions of Theorem \ref{t:uniflogconcBkappa}. Let us introduce some notation. We denote by $S$ the restriction to $\partial K$ of the Lebesgue measure on $\R^d$, and recall $\div_\rho g=\div(\rho g)/\rho$. 
We denote by $n_K$ the outer unit-normal of $K$, and recall that the second fundamental form ${\rm II}$ is defined by 
$$
{\rm II}(X,Y) = \langle \nabla_X n_K,Y\rangle, \qquad X,Y\in T\partial K.
$$
Since $K$ is convex, it satisfies 
${\rm II}(X,X)\geq 0$ for any $X\in T\partial K$. We will make use of the identity
\begin{equation}\label{e:ideSND}
-\langle \nabla_XX, n_K\rangle=\mathrm{II}(X,X)
\end{equation}
on $\partial K$. Its proof goes as follows. Let $x\in\partial K$ and $\gamma$ be a $\calC^1$ curve in $\partial K$ with $\gamma(0)=x$ and $\gamma'(0)=X(x)$. Differentiating the identity $\langle X(\gamma(t)),n_K(\gamma(t))\rangle=0$ at $t=0$ gives
$$
\langle(\nabla_X X)(x),n_K(x)\rangle+\langle X(x),(dn_K)_x(X(x))\rangle=0
$$
which is exactly \eqref{e:ideSND}. 

The proof of Theorem \ref{t:uniflogconcBkappa} is based on the following generalization of Bochner's formula. This generalization is closely linked to Reilly's formula \cite{reilly1977applications}, an extension of Bochner's formula to manifolds with boundary. For more on Reilly's formula, we refer to \cite{kolesnikov2017brascamp} and references therein, where Kolesnikov and Milman generalized further Reilly's formula and provided applications to Brascamp-Lieb type inequalities.
\begin{lemma}[Generalization of the Bochner-Reilly formula] \label{l:genbochner}
Let $g\in \calC^2(K)$ with $\langle g,n_K\rangle=0$ on $\partial K$. Then
\begin{equation}\label{eq:reilly} 
\int (\div_\rho g)^2 \, \dd\rho = 
\int \Tr((\nabla g)^2) \, \dd\rho +\int  \langle g, (\nabla^2 V) g\rangle \dd\rho+\int_{\partial K}\mathrm{II}(g,g)\,e^{-V}\,\dd S .
\end{equation}
\end{lemma}
We formally recover the classical integrated form of Bochner's formula~\eqref{e:bochner} 
from Lemma~\ref{l:genbochner} 
by taking $g$ to be a gradient vector field, and $K=\R^d$ (which is not rigorously allowed in Lemma \ref{l:genbochner} since $K$ is compact, but the proof is the same).
\begin{proof}
We start from the pointwise identity 
\begin{equation}\label{e:pointwisebochner}
\Tr((\nabla g)^2) + \langle g, \nabla^2 V g\rangle - (\div_\rho g)^2= \div_\rho(\nabla_g g - (\div_\rho g)\, g) 
\end{equation}
which is proved for instance in \cite[Lemma~4.1]{antonini2025}. We integrate this identity with respect to $\rho$, use the divergence theorem, the fact that $\langle g,n_K\rangle=0$, and \eqref{e:ideSND}. This yields \eqref{eq:reilly}.
\end{proof}

The generalized Bochner-Reilly formula implies a Poincaré inequality for the (truncated) probability measure $\rho$, with a constant independent of the truncation. Although this Poincaré inequality is well-known, we provide a short proof for clarity of exposition. 
\begin{corollary}[Poincaré inequality] \label{c:poincare}
For any $u\in \dot H^1(\rho)$, there holds $\|\nabla u\|_{L^2(\rho)}^2\geq \kappa \|u\|_{L^2(\rho)}^2$. 
\end{corollary}
\begin{proof}
This follows directly from Lemma \ref{l:genbochner} by taking $g=\nabla h$ with $h$ an eigenfunction corresponding to the first non-trivial eigenvalue of $\Delta_\rho$, as below \eqref{e:bochner}, and using that the boundary term involving the second fundamental form in \eqref{eq:reilly} is non-negative since $K$ is convex.
\end{proof}

\begin{lemma}\label{lem:squaredtrace}
Let $g=(\nabla^2 \phi)^{-1}\nabla u$. Then $\Tr((\nabla g)^2)\geq 0$
over $K$.
\end{lemma}
\begin{proof}
If $g$ were a gradient vector field, the statement would follow from the fact that $\nabla g$ is then a symmetric matrix. In general, it is not the case, but $g \circ \nabla \phi^*$ (the pushforward of $g$ to the target) is a gradient, and this is sufficient to conclude. 
Indeed, let $v=u\circ \nabla \phi^*$. Using \eqref{e:hessphihessphi*}, we get $g=\nabla v\circ \nabla \phi$. Hence 
$\nabla g=(\nabla^2v \circ \nabla \phi) \nabla^2\phi$ and it follows that  
$$
(\nabla^2\phi)^{\frac12}\nabla g(\nabla^2\phi)^{-\frac12}=(\nabla^2\phi)^{\frac12}(\nabla^2v \circ \nabla \phi) (\nabla^2\phi)^{\frac12}
$$
and this last matrix is a symmetric matrix, denoted by $S$. Therefore $\Tr((\nabla g)^2)=\Tr(S^2)\geq 0$.
\end{proof}

\begin{proof}[Proof of Theorem \ref{t:uniflogconcBkappa}] Let  $Y\in \calQ$ and $\mu \in \calF(Y)$. Then, 
one has $\phi_\mu\in \calC^\infty(K)$
according to equation~\eqref{eq:Cinfty_potentials}. 
Let $f\in \calC^2(K)\cap L^2_0(\rho)$, and let $u$ the solution to $E_\mu u =f$. By global Schauder theory for divergence form operators \cite[Chapters 6 and 8, notably Theorem 6.19]{gilbarg2001}, $u\in \calC^{3,\alpha}(K)$ for some $\alpha>0$, and therefore $g = (\nabla^2\phi_\mu)^{-1}\nabla u\in \calC^2(K)$. We apply Lemma \ref{l:genbochner} to $g$ (which verifies $\langle g,n_K\rangle=0$). Using $-\div_\rho g = E_\mu u$, Lemma \ref{lem:squaredtrace},  and the fact that the boundary term in \eqref{eq:reilly} is non-negative (since $K$ is convex), we deduce that 
\begin{equation}\label{e:analogbochner}
\|f\|_{L^2(\rho)}^2 = 
\int (\div_\rho g)^2 \dd\rho 
\geq 
\kappa\|g\|_{L^2(\rho)}^2.
\end{equation}
This is a priori true only for $f\in \calC^2(K)\cap L^2_0(\rho)$, but since the mapping $f\mapsto u$ is continuous from $L^2_0(\rho)$ to $\dot H^1(\rho)$ (see Lemma \ref{lem:weak_solutions}), we get that \eqref{e:analogbochner} holds for any $f\in L^2_0(\rho)$, which concludes the proof of $B(\kappa)$.
\end{proof}

As a consequence, we deduce the following result.
\begin{theorem}\label{t:thwitha}
Let $K\subset\R^d$ be a  uniformly convex body with $\calC^\infty$-boundary, and let $U\in \calC^2(K)$ satisfying $D^2U\succeq \kappa I_d$ and such that 
\begin{equation}\label{e:truncaterho}
\dd\rho=e^{-U}\mathbf1_K\,\dd x
\end{equation}
is a probability measure on $K$.
Then, for any finite convex 
functions $\phi_0,\phi_1$
such that $\nabla\phi_0,\nabla\phi_1 \in L^2(\rho)$, there holds
\begin{equation}\label{e:varkappaW2}
\Var_\rho(\phi_0-\phi_1)
\leq
\frac1\kappa
W_2^2\left((\nabla\phi_0)_\#\rho,(\nabla\phi_1)_\#\rho\right).
\end{equation}
\end{theorem}
\begin{proof}[Proof of Theorem \ref{t:thwitha}]
Combining Theorem \ref{t:uniflogconcBkappa} with Lemma \ref{lem:Ephi}, we obtain that $A(\kappa)$ holds. Since $\rho$ supports a Poincaré inequality by Corollary \ref{c:poincare}, we may then apply Theorem \ref{t:general_reduction} to get Theorem \ref{t:thwitha}.
\end{proof}
Notice that the above variance bound does not involve the convex body $K$. Therefore, we will be able to deduce Theorem \ref{t:main} from Theorem \ref{t:thwitha} by a truncation procedure.
\begin{proof}[Proof of Theorem \ref{t:main}]
If the support of $\rho$ is compact, then Theorem \ref{t:main} immediately follows from Theorem \ref{t:thwitha} together with the regularization arguments described in Appendix \ref{app:regrhoandK}. If the support of $\rho$, denoted by $S$, is not compact, we proceed by truncation. 
For $R>0$, we consider $\rho_R$ the probability measure $\rho$ conditioned on $S_R=B(0,R)\cap S$. Thanks to the regularization arguments of Appendix \ref{app:regrhoandK}, we may furthermore assume that $S_R$ is uniformly convex, with $\calC^\infty$ boundary, and that the restriction $U$ of $V$ to $S_R$ is in $\calC^\infty(S_R)$. Applying Theorem \ref{t:thwitha} to $\rho_R$, with $U$ the restriction of $V$ to $S_R$, we get
\begin{equation}\label{e:brenierstabR}
\Var_{\rho_R}(\phi_0-\phi_1)\leq \frac1\kappa W_2^2\left((\nabla\phi_0)_\#\rho_R,(\nabla\phi_1)_\#\rho_R\right)
\end{equation}
Notice that the constant in the right-hand side does not depend on $R$. 

Let us first prove that 
\begin{equation}\label{e:convW2R}
W_2\left((\nabla\phi_0)_\#\rho_R,(\nabla\phi_1)_\#\rho_R\right)\rightarrow W_2\left((\nabla\phi_0)_\#\rho,(\nabla\phi_1)_\#\rho\right)
\end{equation}
as $R\rightarrow +\infty$. 
Due to the triangle inequality, it is sufficient to prove 
\begin{equation}\label{e:convW2R2}W_2((\nabla\phi_i)_\#\rho_R,(\nabla\phi_i)_\#\rho)\rightarrow 0
\end{equation} 
for $i=0,1$. Let $Z_R=\rho(S_R)$. Then $Z_R\rightarrow 1$ as $R\rightarrow +\infty$. Let $g$ be continuous with $|g(y)|\leq C(1+|y|^2)$. Then
$$
\int g \, \dd (\nabla\phi_i)_\#\rho_R = \frac{1}{Z_R} \int \mathbf{1}_{S_R}g\circ \nabla\phi_i \, \dd\rho \underset{R\rightarrow +\infty}{\longrightarrow} \int g\circ \nabla\phi_i \, \dd\rho = \int g \, \dd(\nabla\phi_i)_\#\rho
$$
by dominated convergence, with dominating function $C(1+|\nabla\phi_i|^2)\in L^1(\rho)$. Since this is true for any $g$ as above, \eqref{e:convW2R2} follows by \cite[Theorem 6.9]{villani2009optimal}.

Let us then prove that
\begin{equation}\label{e:convvarR}
\Var_{\rho_R}(\phi_0-\phi_1)\rightarrow \Var_\rho(\phi_0-\phi_1)
\end{equation}
as $R\rightarrow +\infty$. Letting $f=\phi_0-\phi_1$, we have
$$\Var_{\rho_R}(f)=\frac{1}{2Z_R^2}\iint \mathbf 1_{S_R}(x)\mathbf 1_{S_R}(y)\,(f(x)-f(y))^2\,\dd\rho(x)\,\dd\rho(y).$$
The integrand increases to $(f(x)-f(y))^2$ as $R\rightarrow +\infty$. By monotone convergence and the fact that $Z_R\to1$, we obtain \eqref{e:convvarR}. Notice that $\Var_\rho(\phi_0-\phi_1)<\infty$ as a consequence of the Poincaré inequality \eqref{eq:unif_Cp} for $\rho$, and the fact that $\nabla\phi_0,\nabla\phi_1\in L^2(\rho)$.

Finally, combining \eqref{e:brenierstabR}, \eqref{e:convW2R} and \eqref{e:convvarR}, we get Theorem \ref{t:main}. 
\end{proof}

\subsection{Proof of Corollary \ref{c:bures}} \label{s:bures}

Due to \eqref{e:BWformula}, we can assume without loss of generality that $m_1=m_2=0$. Let $\rho=\mathcal{N}(0,\Id_d)$, and $\phi_i(x)=\frac12 \langle \Sigma_i^{1/2} x,x\rangle$ for $i=1,2$.  First, since $(\nabla \phi_i)_\#\rho = \mathcal{N}(0,\Sigma_i)$, we have by \eqref{e:coupling}
$$
W_2^2(\mathcal{N}(0,\Sigma_1),\mathcal{N}(0,\Sigma_2))\leq \|\nabla\phi_1-\nabla\phi_2\|_{L^2(\rho)}^2= \int |(\Sigma_1^{1/2}-\Sigma_2^{1/2})x|^2\dd\rho(x)=H(\Sigma_1,\Sigma_2)
$$
which gives the right-hand side inequality.
Then we observe that 
$$
\Var_\rho(\phi_1-\phi_2)=\Var_\rho\left(\frac12 \langle (\Sigma_1^{1/2}-\Sigma_2^{1/2}) x,x\rangle\right)=\frac12 H(\Sigma_1,\Sigma_2).
$$
Hence the left-hand side inequality follows from Theorem \ref{t:main}, with $\kappa=1$ since $\rho=\mathcal{N}(0,\Id_d)$.

\section{Proof of Theorem \ref{t:convexset} and Corollary \ref{c:1/3} }\label{s:coercive_cvx}
The proof of Theorem \ref{t:convexset} is based on the following extension of Theorem \ref{t:uniflogconcBkappa} (which corresponds to the case where $a$ is constant, equal to $\kappa^{-1}$).
\begin{theorem}\label{t:thwitha2}
Let $K\subset\R^d$ be a uniformly convex body, with $\calC^\infty$ boundary, and let
$$
\dd\rho=e^{-V(x)}\mathbf1_K\,\dd x
$$
be a probability measure on $K$, where $V\in \calC^2(K)$ is convex.
Suppose there exists $a\in \calC^2(K)$ such that
\begin{equation}\label{e:assumptona}
a\geq0,
\qquad
a\nabla^2V-\nabla^2a\succeq I,
\end{equation}
Let 
$$
\gamma=\left(\sup_K|\nabla a|+\sup_K\sqrt{a+|\nabla a|^2}\right)^{2}
$$
and assume $\gamma<\infty$.
Then, $B(\kappa)$ holds, for $\kappa=\gamma^{-1}$.
\end{theorem}
The function $a$ appearing in this statement helps recovering some coercivity even when $V$ is not strongly convex. For this, one integrates the pointwise formula \eqref{e:pointwisebochner} against $a\, \dd\rho$. A similar proof technique has been used for instance in \cite[Lemma 4.2]{antonini2025} to prove second-order estimates for $p$-Laplace operators with log-concave weights. 

The starting point of the proof of Theorem \ref{t:thwitha2} is the following generalization of Lemma \ref{l:genbochner}.
\begin{lemma} \label{l:genbochner14}
Let $g\in \calC^2(K)$ with $\langle g,n_K\rangle=0$ on $\partial K$. Then
\begin{equation}
\begin{aligned}\label{eq:reilly14} 
\int a(\div_\rho g)^2 \, \dd\rho & = 
\int a\Tr((\nabla g)^2)\,\dd\rho +\int  \langle g, (a\nabla^2 V-\nabla^2 a) g\rangle\, \dd\rho \\
& \qquad - 2 \int (\div_\rho g)\langle \nabla a, g\rangle \, \dd\rho + \int_{\partial K} a\, {\rm II}(g,g) e^{-V} \dd S
\end{aligned}
\end{equation} 
\end{lemma}
\begin{proof}
As in Lemma \ref{l:genbochner}, we start from the pointwise formula \eqref{e:pointwisebochner} and integrate it against $a\,\dd\rho$. The only term that needs manipulations is the one on the right-hand side of \eqref{e:pointwisebochner}.
Using the divergence theorem, $\langle g,n_K\rangle=0$ and \eqref{e:ideSND}, we get
$$
-\int a \, \div_\rho(\nabla_g g - (\div_\rho g)\, g) \, \dd\rho = \int_{\partial K} a\, {\rm II}(g,g) e^{-V} \dd S + \int \langle\nabla a, \nabla_gg-(\div_\rho g)\, g \rangle \dd\rho.
$$
By the product rule, $\langle g,\nabla \langle \nabla a,g\rangle\rangle=\langle \nabla a, \nabla_gg\rangle+\langle g,(\nabla^2a)g\rangle$. Hence
\begin{align*}
\int \langle\nabla a, \nabla_gg\rangle \dd\rho &= \int \langle g,\nabla \langle \nabla a,g\rangle\rangle \dd\rho - \int \langle g,(\nabla^2a)g \rangle \dd\rho \\
&= -\int (\div_\rho g) \langle \nabla a,g\rangle \dd\rho - \int \langle g,(\nabla^2a)g \rangle \dd\rho
\end{align*}
where in the last line we use again the divergence theorem and the fact that $\langle g,n_K\rangle=0$. All in all, we get \eqref{eq:reilly14}.
\end{proof}
 We finally conclude the proof of Theorem \ref{t:thwitha2}.
\begin{proof}[Proof of Theorem \ref{t:thwitha2}]
Let $f\in \calC^{2}(K)\cap L^2_0(\rho)$, and $u$ the solution to $E_\mu u =f$. As in the proof of Theorem \ref{t:uniflogconcBkappa}, $g=(\nabla^2\phi)^{-1}\nabla u \in \calC^2(K)$. We apply Lemma \ref{l:genbochner14} and use \eqref{e:assumptona}, Lemma \ref{lem:squaredtrace}, and the fact that ${\rm II}(g,g)\geq 0$ on $\partial K$ to deduce 
$$
\int |g|^2 \,\dd\rho - 2\int (\div_\rho g)\nabla a\cdot g \, \dd\rho \leq \int a(\div_\rho g)^2\,\dd\rho.
$$
Completing the square and using that $-\div_\rho g=E_\mu u$, we get
$$
\int |(\nabla^2\phi)^{-1}\nabla u+(E_\mu u)\nabla a|^2\, \dd\rho \leq \int (a+|\nabla a|^2)(E_\mu u)^2 \, \dd\rho.
$$
Therefore, by the triangle inequality,
$$
\|(\nabla^2\phi)^{-1}\nabla u\|_{L^2(\rho)} \leq \gamma^{1/2} \|f\|_{L^2(\rho)}.
$$
This is true a priori only for $f\in\calC^2(K)\cap L^2_0(\rho)$, but since $f\mapsto u$ is continuous (see Lemma \ref{lem:weak_solutions}), it holds in fact for any $f\in L^2_0(\rho)$. This concludes the proof of $B(\kappa)$, with $\kappa=\gamma^{-1}$.
\end{proof}

We finally conclude the proofs of Theorem \ref{t:convexset} and Corollary \ref{c:1/3}.
\begin{proof}[Proof of Theorem \ref{t:convexset}]
Let $K\subset B(x_0,R)$ be a convex body, and let
$$
\dd\rho=e^{-V(x)}\mathbf{1}_K(x)\,\dd x,
$$
where $\mathbf{1}_K$ is the characteristic function of $K$, and $V$ is convex. Thanks to the regularization arguments in Appendix \ref{app:regrhoandK}, we may furthermore assume that $K$ has a $\calC^\infty$ boundary, is uniformly convex, and that $V\in\calC^2(K)$. Theorem \ref{t:thwitha2} applied with $a=\frac12(R^2-|x-x_0|^2)$ implies that $B(\kappa)$ holds, with $\kappa=\gamma^{-1}=(4R^2)^{-1}$. According to Lemma \ref{lem:Ephi} this implies that $A(\kappa)$ holds, with same $\kappa$. Moreover, $\rho$ satisfies a Poincaré inequality \cite{payne1960optimal}. Hence, applying Theorem \ref{t:uniflogconcBkappa}, we get Theorem \ref{t:convexset}.
\end{proof}
\begin{proof}[Proof of Corollary \ref{c:1/3}]
Notice that $\phi_{0}$ and $\phi_1$ are both Lipschitz on $K$, with a constant that only depends on $Y$. Applying \cite[Proposition 4.1]{delalande2023quantitative}, we have
$$
\|T_{\mu_0}-T_{\mu_1}\|_{L^2(K,\R^d)}\leq C\|\phi_1-\phi_0\|_{L^2(K)}^{1/3}
$$
where the $L^2$ norms are with respect to the Lebesgue measure on $K$. Since the density of $\rho$ is upper and lower bounded on $K$, this yields
$$
\|T_{\mu_0}-T_{\mu_1}\|_{L^2(\rho,\R^d)}\leq C\|\phi_1-\phi_0\|_{L^2(\rho)}^{1/3}
$$
and we conclude the proof thanks to Theorem \ref{t:convexset}.
\end{proof}

\section{Proof of Theorem \ref{t:directionalpoincare} and Corollary \ref{c:failkls}}
\label{sec:negative}
In this section, we prove Theorem \ref{t:directionalpoincare} and Corollary \ref{c:failkls}, which show in particular that not all log-concave probability
measures 
satisfy a Brenier stability inequality.
For $x\in\R^d$, we write $x=(x_1,\ldots,x_d)$. 

\begin{proof}[Proof of Theorem \ref{t:directionalpoincare}]
Without loss of generality, let us establish the result for $j=d$. Given $f\in C_c^\infty(\R^d)$, odd in the $d$-th variable, and $\delta>0$, we consider 
$$
\phi_{\delta}(x)=\frac12(x_1^2+\ldots+x_{d-1}^2+\delta x_d^2), \qquad \phi_{\delta,\varepsilon}=\phi_\delta+\varepsilon f
$$
which both have gradients
in $L^2(\rho)$.
For each fixed $f,\delta$, the second potential $\phi_{\delta,\varepsilon}$ is uniformly convex for sufficiently small $|\varepsilon|$. Let us also consider a compactly supported vector field $b$  such that $\div(\rho b)=0$. Its components, denoted by $b_1,\ldots,b_d$,  will be chosen later in terms of $f$ and $\rho$. Let $(F_t)_{t\in\R}$ be the flow of $b$. It verifies $(F_t)_\#\rho = \rho$ for any $t\in\R$, therefore $(\nabla\phi_\delta\circ F_t,\nabla \phi_{\delta,\varepsilon})_\#\rho$ is a coupling of $(\nabla\phi_\delta)_\#\rho$ and of $(\nabla\phi_{\delta,\varepsilon})_\#\rho$. Hence, 
$$
W_2^2((\nabla\phi_\delta)_\#\rho,(\nabla\phi_{\delta,\varepsilon})_\#\rho)\leq \|\nabla\phi_\delta\circ F_t - \nabla\phi_{\delta,\varepsilon}\|_{L^2(\rho)}^2.
$$
We choose $t=\varepsilon$ and let $\varepsilon\rightarrow 0$. The components of the vector field $\nabla\phi_\delta\circ F_\varepsilon$ are
$$
(\nabla\phi_\delta\circ F_\varepsilon)_i(x)= \begin{cases} x_i+\varepsilon b_i+o(\varepsilon)  & \text{for } \ i=1,\ldots,d-1 \\  \delta (x_d+\varepsilon b_d+o(\varepsilon)) & \text{for } \ i=d. \end{cases}
$$
Hence
\begin{equation}\label{e:limsupW2}
\limsup_{\varepsilon\rightarrow 0} \frac{W_2^2((\nabla\phi_\delta)_\#\rho,(\nabla\phi_{\delta,\varepsilon})_\#\rho)}{\varepsilon^2}\leq \int \Bigl( |\delta b_d-\partial_d f|^2 + \sum_{i=1}^{d-1} |b_i-\partial_i f|^2 \Bigr) \dd\rho.
\end{equation}
We choose $b_i=\partial_if$ for $i=1,\ldots,d-1$ to cancel the sum, and $b_d$ in a way that $\div(b\rho)=0$. In other words, identifying the probability measure $\rho$ with its density, we let
$$
b_d(x)=-\frac{1}{\rho(x)}\int_{-\infty}^{x_d}\sum_{i=1}^{d-1} \partial_i\big(\rho(x_1,\ldots,x_{d-1},z)\partial_i f(x_1,\ldots,x_{d-1},z)\big)
\dd z.
$$
Since $f$ is odd in the variable $x_d$ and $\rho$ is even in this variable (because it is unconditional), we get that the vector field $b$ is compactly supported: indeed, the integrand in the definition of $b_d$ is odd in $z$, thus the integral over the whole line is equal to $0$, and $b$ vanishes outside the support of $f$.
Sending $\delta\rightarrow 0$ in \eqref{e:limsupW2} we obtain 
$$
\limsup_{\varepsilon\rightarrow 0} \frac{W_2^2((\nabla\phi_\delta)_\#\rho,(\nabla\phi_{\delta,\varepsilon})_\#\rho)}{\varepsilon^2}\leq \int |\partial_d f|^2 \dd\rho. 
$$ 
Since $\Var_\rho(\phi_\delta-\phi_{\delta,\varepsilon})=\varepsilon^2\Var_\rho(f)$, this concludes the proof of the result.
\end{proof}

\begin{proof}[Proof of Corollary \ref{c:failkls}]
Consider the unconditional and log-concave probability density on $\R^d$ given by
\begin{equation}\label{e:rhocontreexemple}
\dd\rho(x)=\frac1Z e^{-\sqrt{1+|x|^2}}\dd x
\end{equation}
where $Z$ is a normalizing constant.  Observe that $\sqrt{1+|x|^2}$ is strictly convex but not uniformly  convex since it has eigenvalues approaching $0$ along the axis $x_2=0$.  Let us prove that the inequality \eqref{eq:directional_poincare} fails for any $C_W(\rho)>0$ (and for instance $j=d$), which will prove Corollary \ref{c:failkls} according to Theorem \ref{t:directionalpoincare}. The argument described below works for many other log-concave probability densities than \eqref{e:rhocontreexemple}, for instance $\dd\rho(x)=\frac1Z e^{-|x|}\dd x$.

For $x=(x_1,\ldots,x_d)\in\R^d$, we write $x'=(x_1,\ldots,x_{d-1})$, hence $x=(x',x_d)$. We observe that on sets of the form 
$$
\Box_R=\left\{x=(x',x_d)\in\R^d \mid R\leq |x'|\leq R+1,\ |x_d|\leq \sqrt{R}\right\}
$$
for large $R>0$, the density of $\rho$ does not vary much. We show that  \eqref{eq:directional_poincare} fails using the fact that the Poincaré constant of the uniform measure on $\Box_R$ tends to $0$ as $R\rightarrow +\infty$.
More precisely, let $f\in C_c^\infty((0,1))$ and $g\in C_c^\infty((-1,1))$ with $g$ odd and define
$$
h_R(x)=f(|x'|-R)g(x_d/\sqrt{R})
$$
which is supported in $\Box_R$ and has mean $0$. For any $x=(x',x_d)\in \Box_R$, witing $x'=(R+s)\frac{x'}{|x'|}$ and $x_d=t\sqrt{R}$ where $s\in[0,1]$ and $t\in [-1,1]$, we have
$$
\sqrt{1+|x|^2}=\sqrt{1+(R+s)^2+Rt^2}=R+s+\frac{t^2}{2}+O(R^{-1}).
$$
Hence, as $R\rightarrow +\infty$,  
$$
\int h_R^2\, \dd\rho \sim Z'^{-1}e^{-R}\sqrt{R}\left(\int f(s)^2s^{d-2}e^{-s}\dd s\right)\left(\int g(t)^2e^{-\frac{t^2}{2}}\dd t\right)
$$
and
$$
\int |\partial_d h_R|^2 \dd\rho = Z'^{-1}\frac{e^{-R}}{\sqrt{R}}\left(\int f(s)^2s^{d-2}e^{-s}\dd s\right)\left(\int g'(t)^2e^{-\frac{t^2}{2}}\dd t\right).
$$
As a consequence, for any $C>0$, 
$$
\Var_\rho(f)\leq C_W(\rho) \int \|\partial_d f\|_{L^2(\rho)}^2
$$
cannot hold. 
\end{proof}

\appendix

\section{Omitted Proofs}
\label{a:caracW2}

\subsection{Proof of the identity~\eqref{e:W2asinf}} 
In this brief appendix, we provide
a self-contained proof
of the identity~\eqref{e:W2asinf}. 
\begin{proposition}\label{p:W2asinf}
Let $\rho\in\calP_2(\bbR^d)$
be absolutely continuous 
with respect to the Lebesgue
measure on $\bbR^d$. Then,
for all convex functions $\phi_0,\phi_1$
such that $\nabla\phi_0,\nabla\phi_1\in L^2(\rho)$,
$$W_2((\nabla\phi_0)_\#\rho,(\nabla\phi_1)_\#\rho) = \inf_{S \in \calS(\rho)} 
\|\nabla\phi_1 - \nabla\phi_0\circ S\|_{L^2(\rho)}.$$
\end{proposition}
\begin{proof} 
Define 
$$\mu_0= (\nabla\phi_0)_\#\rho,
\quad \mu_1 = (\nabla\phi_1)_\# \rho.$$
Let us begin by assuming
that the functions $\phi_0,\phi_1$
are strongly convex, in which case the
measures $\mu_0$ and $\mu_1$
are absolutely continuous with respect
to the Lebesgue measure on $\bbR^d$, and $\phi_i^*$ is finite with a Lipschitz gradient.
 
Consider an optimal transport map $T$ from $\mu_1$ to $\mu_0$, which exists according to Brenier's theorem. Set $S=\nabla \phi_0^*\circ T\circ \nabla\phi_1$. This map satisfies $\rho$-a.e the identity $\nabla\phi_0\circ S = T\circ\nabla\phi_1$, and $S\in\mathcal{S}(\rho)$. Therefore
$$
\|\nabla\phi_1 - \nabla\phi_0\circ S\|_{L^2(\rho)}=\|\nabla\phi_1-T\circ\nabla\phi_1\|_{L^2(\rho)}=\|{\rm id}-T\|_{L^2(\nu)}=W_2(\mu,\nu).
$$
Now, for any $S\in\mathcal{S}(\rho)$, the coupling $(\nabla\phi_1,\nabla\phi_0\circ S)_\#\rho$ is a coupling between $\nu$ and $\mu$, with cost $\|\nabla\phi_1-\nabla\phi_0\circ S\|_{L^2(\rho)}$ and therefore 
$$
W_2(\mu_0,\mu_1)\leq \|\nabla\phi_1-\nabla\phi_0\circ S\|_{L^2(\rho)}.
$$
This proves the claim when $\phi_0$ and $\phi_1$
are strongly convex. 

To deduce the general claim, let 
$\phi_0,\phi_1$ be arbitrary
convex functions with gradients
in $L^2(\rho)$, and let $\phi_i^\varepsilon = \phi_i + \frac{\varepsilon}{2}|x|^2$ for $\varepsilon>0$, which is $\varepsilon$-strongly convex. Let us pass to the limit $\varepsilon\rightarrow 0$ in the formula
$$ W_2((\nabla\phi_0^\varepsilon)_\#\rho,(\nabla\phi_1^\varepsilon)_\#\rho) = \inf_{S \in \calS(\rho)} 
\|\nabla\phi_1^\varepsilon - \nabla\phi_0^\varepsilon\circ S\|_{L^2(\rho)}.$$
Let $\mu_i^\varepsilon := (\nabla\phi_i^\varepsilon)_\#\rho$. We have
$$
W_2(\mu_i^\varepsilon,\mu_i)\leq \|\nabla\phi_i^\varepsilon-\nabla\phi_i\|_{L^2(\rho)}= \varepsilon \left(\int |x|^2 \,\dd\rho(x)\right)^{1/2}
$$
hence $W_2(\mu_0^\varepsilon,\mu_1^\varepsilon)\rightarrow W_2(\mu_0,\mu_1)$ as $\varepsilon\rightarrow 0$. For $S\in\calS(\rho)$, we have 
$$
\|\nabla\phi_0^\varepsilon \circ S-\nabla\phi_0\circ S\|_{L^2(\rho)}=\|\nabla\phi_0^\varepsilon-\nabla\phi_0\|_{L^2(\rho)}
$$
hence $\|\nabla\phi_1^\varepsilon - \nabla\phi_0^\varepsilon\circ S\|_{L^2(\rho)}\rightarrow \|\nabla\phi_1 - \nabla\phi_0\circ S\|_{L^2(\rho)}$ as $\varepsilon\rightarrow 0$. This proves Proposition \ref{p:W2asinf} in full generality.
\end{proof}

\subsection{Regularization of $\rho$ and $K$} \label{app:regrhoandK}
In the proofs of both Theorem \ref{t:main} and Theorem \ref{t:convexset}, we assumed that $V$ is $\calC^2$, and in the proof of Theorem \ref{t:convexset} we also assumed that $K$ has $\calC^\infty$ boundary and is uniformly convex. In this appendix, we explain how to reduce to this case. We focus on the case of Theorem \ref{t:convexset} since the reduction for Theorem \ref{t:main} is easier and follows from a straightforward adaptation. 

In the sequel, we assume Theorem \ref{t:convexset} has been proved under the additional assumptions that $V\in \calC^2(K)$ and $K$ has $\calC^\infty$ boundary and is uniformly convex, and we show this implies Theorem \ref{t:convexset} in full generality. Without loss of generality, $V$ is finite in the interior of $K$, which we denote by $\Omega$.
Up to a translation, we may assume $B(0,r)\subset\Omega$ for some $r>0$. We let 
$$
p(x) = \inf\{t>0 : x\in t\Omega\}
$$
be the gauge of $\Omega$.
It is convex, $L$-Lipschitz with $L = 1/r$, and $\Omega = \{p<1\}$.
Consider two convex functions $\phi_0$ and $\phi_1$ such that $\phi_0,\phi_1\in \calC^{1,\alpha}(K)$ for some $\alpha>0$. 

Let $\chi\in C_c^\infty(B(0,1))$, and $\chi_\delta=\delta^{-d}\chi(\delta^{-1}\cdot)$ for $\delta>0$. Let $p_\delta = p*\chi_\delta+\delta|x|^2$, which is smooth and convex, with $|p_\delta - p|\le \delta(L+\sup_{x\in K} |x|^2)$ in $K$. Set
$$K_\delta = \{p_\delta \le 1 - 2L\delta\},$$
which is compact, convex, and included in $\Omega$.
For $0<\delta<\delta_0$ with $\delta_0$ small enough, we have $p_\delta> \min p_\delta$ on $\partial K_\delta$, so no boundary point minimizes the convex function $p_\delta$. Hence $\nabla p_\delta \neq 0$ on $\partial K_\delta$, and $\partial K_\delta$ is $\calC^\infty$ by the implicit function theorem. Moreover, $K_\delta$ is uniformly convex, as a sublevel set of the strongly convex function $p_\delta$. 

Let $\delta\in (0,\delta_0)$ be fixed for a moment. Let $\rho_\delta$ be the probability measure whose density is proportional to $\mathbf{1}_{K_\delta}e^{-V}$. For $\varepsilon<{\rm dist}(K_\delta,\partial \Omega)$, let $V_\varepsilon=V\star \chi_\varepsilon$, which is $\calC^\infty$ in a neighborhood of $K_\delta$, and convex. Applying Theorem \ref{t:convexset} to the probability measure $\rho_{\delta,\varepsilon}$ whose density is proportional to $\mathbf{1}_{K_\delta}e^{-V_\varepsilon}$, we get
$$
\Var_{\rho_{\delta,\varepsilon}}(\phi_1-\phi_0)\leq 4R^2W_2^2((\nabla\phi_0)_\#\rho_{\delta,\varepsilon},(\nabla\phi_1)_\#\rho_{\delta,\varepsilon}).
$$
Let us show that taking the limit $\varepsilon\rightarrow 0$, this implies 
\begin{equation}\label{e:avecdelta}
\Var_{\rho_{\delta}}(\phi_1-\phi_0)\leq 4R^2W_2^2((\nabla\phi_0)_\#\rho_{\delta},(\nabla\phi_1)_\#\rho_{\delta}).
\end{equation}
Since $V_\varepsilon \to V$ uniformly on $K_\delta$, and $V$ is bounded on $K_\delta$, the density $\rho_{\delta,\varepsilon}$ converges uniformly on $K_\delta$  to $\rho_\delta$ as $\varepsilon\rightarrow 0$, and stays bounded above and below. Hence in the limit $\varepsilon\rightarrow 0$ we obtain $\Var_{\rho_{\delta,\varepsilon}}(\phi_1-\phi_0)\rightarrow \Var_{\rho_\delta}(\phi_1-\phi_0)$. Now we observe that $W_2^2((\nabla\phi_i)_\#\rho_{\delta,\varepsilon},(\nabla\phi_i)_\#\rho_{\delta})\rightarrow 0$ as $\varepsilon\rightarrow 0$ since $\phi_i\in \calC^{1,\alpha}(K)$ (see for instance \cite[Proposition in Section 1.3.1]{carlier2025quantitative}), and \eqref{e:avecdelta} follows. For the same reason, we also get 
$$
W_2^2((\nabla\phi_0)_\#\rho_{\delta},(\nabla\phi_1)_\#\rho_{\delta})\rightarrow W_2^2((\nabla\phi_0)_\#\rho,(\nabla\phi_1)_\#\rho)
$$
as $\delta\rightarrow 0$. Since $\Var_{\rho_\delta}(\phi_1-\phi_0)\rightarrow \Var_\rho(\phi_1-\phi_0)$ as $\delta\rightarrow 0$ by dominated convergence, we finally get
\begin{equation}\label{e:onyestpresque}
\Var_{\rho}(\phi_1-\phi_0)\leq 4R^2W_2^2((\nabla\phi_0)_\#\rho,(\nabla\phi_1)_\#\rho)
\end{equation}
by letting $\delta\rightarrow 0$ in \eqref{e:avecdelta}. We have obtained \eqref{e:onyestpresque} for any convex $\phi_0,\phi_1\in \calC^{1,\alpha}(K)$. Finally, this extends to any convex functions $\phi_0,\phi_1$ 
with gradients in $L^2(\rho)$ by the same arguments as in the proof of Theorem \ref{t:general_reduction}: first reduce to the case where $(\nabla\phi_0)_\#\rho$ and $(\nabla\phi_1)_\#\rho$ have smooth positive densities on their support, which is a ball, then extend to the case where $(\nabla\phi_0)_\#\rho$ and $(\nabla\phi_1)_\#\rho$ are only compactly supported, and then to any $\phi_0,\phi_1$
with $L^2(\rho)$ gradients.

\bibliographystyle{plain}
\bibliography{references}

\end{document}